\documentclass[a4paper,11pt,reqno]{amsart}
\usepackage[centering,includeheadfoot,margin=2.5 cm]{geometry}
\usepackage{graphics,enumitem,epsfig,textcomp}
\usepackage{amsfonts,amsmath,amssymb,euscript,color,mathrsfs}
\usepackage[utf8]{inputenc}	

\usepackage[symbol]{footmisc}

\usepackage
[
pdfauthor   = {Lisa Maria Kreusser},
pdftitle    = {Title},
pdfsubject  = {},
pdfkeywords = {},
bookmarks=true,
bookmarksopen=true,
pagebackref=false,
hyperindex=true,
colorlinks=true,
linkcolor=blue,
citecolor=blue,
filecolor=blue,
urlcolor=blue,
]
{hyperref}					
\usepackage{textcomp}	
\usepackage{color}			
\usepackage{graphicx}		

\usepackage[caption=false,subrefformat=parens,labelformat=parens]{subfig}

\newcommand\N{\mathbb{N}}
\newcommand\Z{\mathbb{Z}}
\newcommand\R{\mathbb{R}}

\newcommand\E{\mathcal{E}}
\newcommand\bl{\left(}
\newcommand\br{\right)}

\newcommand\domain{\Omega}

\newcommand*\di{\mathop{}\!\mathrm{d}}

\newcommand\supp{\operatorname{supp}}

\renewcommand\epsilon{\varepsilon}
\renewcommand\theta{\vartheta}

\newtheoremstyle{mytheoremstyle}
{6pt}                   
{6pt}                    
{\itshape}                  
{}						
{\bf}              
{}                     
{.5em}                    
{}  
\theoremstyle{mytheoremstyle}
\newtheorem{satz}{Satz}[section]
\newtheorem{lemma}[satz]{Lemma}
\newtheorem{proposition}[satz]{Proposition}
\newtheorem{corollary}[satz]{Corollary}
\newtheorem{theorem}[satz]{Theorem}
\newtheorem{assumption}[satz]{Assumption}
\newtheorem{remark}[satz]{Remark}

\newtheoremstyle{mytdefintionstyle} 
{6pt}                
{6pt}                 
{\rm}                 
{}				
{\bf}                 
{}                        
{.5em}                    
{}  
\theoremstyle{mytdefintionstyle}

\graphicspath{{./Numerics/}}

\numberwithin{equation}{section}
\usepackage{cleveref}

\newtheorem{assump}{}

\newcommand\assref[1]{\textbf{\ref{#1}}}

\title{Equilibria of an anisotropic nonlocal interaction equation: Analysis and numerics}

\begin{document}
	
	\maketitle
	
	\centerline{
		{\large Jos\'{e} A. Carrillo}\footnote{Mathematical Institute, University of Oxford, Andrew Wiles Building, Radcliffe Observatory Quarter, Woodstock Road, Oxford OX2 6GG, United Kingdom; 
			{\it carrillo@maths.ox.ac.uk}}\quad
		{\large Bertram D\"uring}\footnote{Mathematics Institute, University of Warwick, Zeeman Building, Coventry
			CV4 7AL, United Kingdom;
			{\it bertram.during@warwick.ac.uk}}\quad
		{\large Lisa Maria Kreusser}\footnote{Department of Applied Mathematics and Theoretical Physics (DAMTP), University of Cambridge, Wilberforce Road, Cambridge CB3 0WA, United Kingdom;
			{\it L.M.Kreusser@damtp.cam.ac.uk}}\quad
		{\large Carola-Bibiane Sch\"{o}nlieb}\footnote{Department of Applied Mathematics and Theoretical Physics (DAMTP), University of Cambridge, Wilberforce Road, Cambridge CB3 0WA, United Kingdom; {\it C.B.Schoenlieb@damtp.cam.ac.uk}}	
	}
	\vskip 10mm

	\noindent{\bf Abstract.}
	In this paper, we study the equilibria of an anisotropic, nonlocal aggregation equation with nonlinear diffusion which does not possess a gradient flow structure. Here, the anisotropy is induced by an underlying tensor field. Anisotropic forces cannot be associated with a potential in general and stationary solutions of anisotropic aggregation equations generally cannot be regarded as minimizers of an energy functional. We derive equilibrium conditions for stationary line patterns in the setting of spatially homogeneous tensor fields. The stationary solutions can be regarded as the minimizers of a regularised energy functional  depending on a scalar potential. A dimension reduction from the two- to the one-dimensional setting allows us to study the associated one-dimensional problem instead of the two-dimensional setting.  We establish $\Gamma$-convergence of the regularised energy functionals as the diffusion coefficient vanishes, and prove the convergence of minimisers of the regularised energy functional to minimisers of the non-regularised energy functional. Further, we investigate properties of stationary solutions on the torus, based on known results in  one spatial dimension. Finally, we prove weak convergence of a numerical scheme for the numerical solution of the anisotropic, nonlocal aggregation equation with nonlinear diffusion and any underlying tensor field, and show numerical results.
	
	\vskip 3mm

	\vskip 7mm

	\section{Introduction}\label{sec:intro}	
	
	The derivation, analysis and numerics of mathematical models for collective behaviour of cells, animals or humans have recently been receiving increasing attention. Based on agent-based modelling approaches,  a variety of continuum models has been derived and used to describe biological aggregations such as flocks and swarms \cite{nonlocal_swarm,TBL}.
	Motivated by the simulation of fingerprint patterns which can be modelled as the interaction of a large number of cells \cite{patternformationanisotropicmodel,Merkel}, a continuum model can be derived following the procedure in \cite{meanfieldlimit,Golse}. The continuum model \cite{patternformationanisotropicmodel} is given by the anisotropic aggregation equation
	\begin{align}\label{eq:macroscopiceq}
	\begin{split}
	\partial_t \rho(t,(x,y))+\nabla\cdot \left[ \rho(t,(x,y))(F(\cdot,T(x,y)) \ast \rho(t,\cdot))(x,y)\right]=0\qquad \text{in }\R_+\times \R^2
	\end{split}
	\end{align} 
	with  initial condition $\rho|_{t=0}  = \rho^{in}$ in $\R^2$ for some given initial data $\rho^{in}$. Here, 
	\begin{align}\label{eq:velocityfield}
	\begin{split}
	u_\rho(t,(x,y))&=(  F(\cdot,T(x,y))\ast \rho(t,\cdot))(x,y)\\& =\int_{\R^2} F((x-w,y-z),T(x,y)) \rho(t,(w,z))\di (w,z)
	\end{split}
	\end{align}	
	is the velocity field with $|u_\rho(t,(x,y))| \leq f$ for the uniform bound $f$ of $F$ where the term  $F((x-w,y-z),T(x,y))$  denotes the force which a particle at position $(w,z)$ exerts on a particle at position $(x,y)$. The left-hand side of \eqref{eq:macroscopiceq} represents the active transport of the density $\rho$ associated to a nonlocal velocity field $u_\rho$.
	
	The force $F$  depends on an underlying stress tensor field $T(x,y)$ at location $(x,y)$.
	The existence of such a tensor field $T(x,y)$  is  motivated by experimental results for  simulating  fingerprints  \cite{KIM1995411} and a model describing the formation of fingerprint patterns based on the interaction of so-called Merkel cells has been suggested by K\"ucken and Champod \cite{Merkel}. In the following, we make general assumptions on the forces which include the explicit choice of forces  suggested in \cite{Merkel}. This more general definition of the forces can be regarded as the starting point for understanding anisotropic pattern formation in nature.
	Since an alignment of mass along the local stress lines is observed, we define the tensor field $T(x,y)$ by the  directions of smallest stress  at location $(x,y)$, i.e.\ we consider a unit vector field  $s=s(x,y)\in\R^2$ and introduce a corresponding orthonormal vector field $l=l(x,y)\in\R^2$, representing the directions of largest stress. The tensor field $T(x,y)$ at $(x,y)$ is given by
	\begin{align}\label{eq:tensorfield}
	T(x,y):=\chi s(x,y)\otimes s(x,y) +l(x,y)\otimes l(x,y)\in\R^{2,2}.
	\end{align}
	The parameter $\chi\in[0,1]$ in the definition of the tensor field introduces an anisotropy in the direction $s$. 
	
	A typical aspect of aggregation models is the competition of social interactions (repulsion and attraction) between the particles which is also the focus of our research. Hence, we assume that the total force $F$  is given by
	\begin{align}\label{eq:totalforce}
	F((x-w,y-z),T(x,y))=F_A((x-w,y-z),T(x,y))+F_R(x-y,w-z).
	\end{align} 
	Here, $F_R$ denotes the repulsion force that a particle at location $(w,z)$ exerts on particle at location $x$ and $F_A$ is the attraction force a particle at location $(w,z)$ exerts on particle at location $(x,y)$. The repulsion and attraction forces are of the form
	\begin{align*}
	F_R(x-w,y-z)=f_R(\sqrt{(x-w)^2+(y-z)^2})\begin{pmatrix}x-w\\ y-z\end{pmatrix}
	\end{align*}
	and
	\begin{align*}
	F_A((x-w,y-z),T(x,y))=f_A(\sqrt{(x-w)^2+(y-z)^2})T(x,y)\begin{pmatrix}x-w\\ y-z\end{pmatrix},
	\end{align*}
	respectively, with radially symmetric coefficient functions $f_R$ and $f_A$, where $(x,y)$, $(w,z)\in\R^2$. An example for the force coefficients $f_R$ and $f_A$ was suggested by
	K\"ucken and Champod \cite{Merkel}, given by
	\begin{align}\label{eq:repulsionforcemodel}
	f_R(\tau)=(\alpha \tau^2+\beta)\exp(-e_R \tau)
	\end{align}
	and 
	\begin{align}\label{eq:attractionforcemodel}
	f_A(\tau)=-\gamma\tau \exp(-e_A\tau)
	\end{align}
	for nonnegative constants $\alpha$, $\beta$, $\gamma$, $e_A$ and $e_R$, and $\tau \geq 0$.
	We assume that the total force \eqref{eq:totalforce} exhibits short-range repulsion and long-range attraction along $l$, and only repulsion along $s$, while the direction of the interaction forces is determined by the parameter $\chi\in[0,1]$ in the definition of  $T$ in \eqref{eq:tensorfield}. These assumptions on the force coefficients are satisfied for the parameters proposed in \cite{During2017}, given by
	\begin{align}\label{eq:parametervaluesRepulsionAttraction}
	\begin{split}
	\alpha&=270, \quad \beta=0.1, \quad \gamma=10.5, \quad
	e_A=95, \quad e_R=100, \quad \chi= 0.2.
	\end{split}
	\end{align} 
	Motivated by plugging \eqref{eq:tensorfield} into the definition of the total force \eqref{eq:totalforce},  we consider a more general form of the total force, given by
	\begin{align}\label{eq:totalforcenew}
	\begin{split}
	F((x-w,y-z),T(x,y))&=f_s(\sqrt{(x-w)^2+(y-z)^2})\left(s(x,y)\cdot \begin{pmatrix}x-w\\y-z\end{pmatrix}\right)s(x,y)\\&\quad+f_l(\sqrt{(x-w)^2+(y-z)^2})\left(l(x,y)\cdot \begin{pmatrix}x-w\\y-z\end{pmatrix}\right)l(x,y)
	\end{split}
	\end{align}
	for coefficient functions $f_s$ and $f_l$, where $f_s=f_R+\chi f_A$  and $f_l=f_R+ f_A$ for the K\"ucken-Champod model. 
	
	The macroscopic model \eqref{eq:macroscopiceq} can be regarded as the rigorous macroscopic limit of an anisotropic particle model	as the number of particles $N$ goes to infinity. The $N$ interacting particles with positions $x_j=x_j(t)\in\R^2$, ${j=1,\ldots,N},$ at time $t$   satisfy
	\begin{align}\label{eq:particlemodel}
	\frac{\di x_j}{\di t}=\frac{1}{N}\sum_{\substack{k=1\\k\neq j}}^N F(x_j-x_k,T(x_j)),
	\end{align}
	equipped with initial data $x_j(0)=x_j^{in},~j=1,\ldots,N$,  for given scalars $x_j^{in},~j=1,\ldots,N$.  A special instance of this model  has been  introduced in \cite{Merkel} for simulating fingerprint patterns where the authors assumed a specific form of the force $F$. The particle model in its general form \eqref{eq:particlemodel} has been studied in \cite{patternformationanisotropicmodel,stabilityanalysisanisotropicmodel,During2017}. The existence of different kinds of steady states, including steady states in the form of lines, is investigated in \cite{patternformationanisotropicmodel}, both for the particle model \eqref{eq:particlemodel} and its continuum counterpart \eqref{eq:macroscopiceq}. The stationary solutions to \eqref{eq:macroscopiceq} can be regarded as solutions with one-dimensional support  \cite{patternformationanisotropicmodel} and may be constant on its support. The direction of the line patterns depends on the choice of the tensor field $T$ with its vector fields $s$ and $l$. For purely repulsive forces along $s$ and  short-range repulsive, long-range attractive forces along $l$, the stability of line patterns is proven for spatially homogeneous tensor fields in \cite{stabilityanalysisanisotropicmodel}, based on a stability analysis of \eqref{eq:particlemodel}. The proof considers perturbations of  equidistantly distributed particles along lines and shows that line patterns along $s$ are stable, while most other rotations including line patterns along $l$ are unstable. This motivates to study constant stationary solutions along $s$ for stable stationary solutions of \eqref{eq:macroscopiceq}.
	The numerical simulations of \eqref{eq:particlemodel} for spatially inhomogeneous tensor fields demonstrate that line patterns can be obtained as stationary solutions \cite{During2017}, again aligned along $s=s(x)$. Applications of \eqref{eq:particlemodel} include the simulation of fingerprints where $s$ is regarded as an underlying stress field. 
	
	Since our fingerprint lines do not have a one-dimensional support and, in fact, have a certain width, we modify \eqref{eq:macroscopiceq}, studied in \cite{patternformationanisotropicmodel,stabilityanalysisanisotropicmodel}. We introduce a small nonlinear diffusion on the right-hand side of \eqref{eq:macroscopiceq} to widen the support  of the line structures. This leads to the nonlocal aggregation equation with nonlinear diffusion
	\begin{align}\label{eq:macroscopiceqnonlin}
	\begin{split}
	&\partial_t \rho(t,(x,y))+\nabla\cdot \left[ \rho(t,(x,y))(F(\cdot,T(x,y)) \ast \rho(t,\cdot))(x,y)\right]\\&=\delta \nabla \cdot (\rho(t,(x,y)) \nabla \rho(t,(x,y))) 
	\end{split}
	\end{align} 
	where $\delta\ll 1$. 
	In particular, for the spatially homogeneous tensor field $T$ with $s=(0,1)$ and $l=(1,0)$  straight vertical lines are obtained as stationary solutions \cite{patternformationanisotropicmodel,stabilityanalysisanisotropicmodel,During2017} which can be regarded as constant solutions  along the vertical axis. For solutions of this form, the diffusion term only acts perpendicular to the  line patterns and not parallel.
	Hence, a positive diffusion coefficient $\delta$ leads to nonlinear diffusion along the horizontal axis and  we expect the  widening of the vertical line profile.  
	
	\subsection{Isotropic aggregation equations}
	While we consider anisotropic aggregation equations of the form \eqref{eq:macroscopiceq} in this work, mainly isotropic aggregation equations  of the form
	\begin{align}\label{eq:standardmodelmacroscopic}
	\rho_t+\nabla\cdot (\rho (-\nabla  W\ast \rho))=0 \qquad \text{in }\R_+\times \R^d
	\end{align}
	have been studied in the literature. Here, $W$ is a radially symmetric interaction potential satisfying  $F=-\nabla W$ on $\R^d$. The study of the isotropic aggregation equations in  terms of its gradient flow structure, the blow-up dynamics for fully attractive potentials, and the rich variety of steady states has attracted the interest of many research groups recently. In these works, the energy
	\begin{align}\label{eq:energyisotropic}
	\E(\rho)= \frac{1}{2}\int_{\R^d} \int_{\R^d} W(u-v) \di \rho(u)\di \rho(v)
	\end{align} 
	in the $d$-dimensional setting plays an important role since it governs the dynamics, and  its (local) minima describe the long-time asymptotics of solutions. Sharp conditions for the existence of global minimisers for a broad class of nonlocal interaction energies on the space of probability measures have been  established in \cite{Simione2015}.
	
	In terms of biological applications, nonlocal interactions  on different scales  are considered for describing the interplay between  short-range repulsion which prevents collisions between  individuals, and long-range attraction  which keeps the swarm cohesive. These repulsive-attractive potentials  can be considered as a minimal model for pattern formation in large systems of individuals \cite{nonlocalinteraction}.	
	
	Very few numerical schemes apart from particle methods have been proposed to simulate solutions of isotropic aggregation equations after blow-up. The so-called sticky particle method \cite{carrillo2011} is a convergent  numerical scheme, used to obtain qualitative properties of the solution such as the finite time total collapse. While numerical results have been obtained in the one-dimensional setting \cite{James2013}, this method is not  practical to deal with finite time blow-up and the behavior of solutions after blow-up in dimensions larger than one. Let the solution to \eqref{eq:macroscopiceq} with initial data $\rho^{in}$ be denoted by $\rho$ and the solution of the particle model \eqref{eq:particlemodel} with initial data $\rho^{in,N}$ be denoted  by $\rho^N(t)=\frac{1}{N}\sum_{j=1}^N \delta(u-x_j(t))$  at time $t\geq 0$ and location $u\in\R^d$. If $F=-\nabla W$ for some radially symmetric potential $W$ and the initial data satisfies $d_W(\rho^{in},\rho^{in,N})\to 0$ as $N\to \infty$ in the Wasserstein distance $d_W$, then
	\begin{align*}
	\sup_{t\in[0,T]} d_W(\rho(t),\rho^N(t)) \to 0.	
	\end{align*}
	for any given $T>0$ \cite{Carrillo2016304}. From the theoretical viewpoint, this is a very nice result, but  in practice a very large number of particles is required for numerical simulations of the particle model \eqref{eq:particlemodel} to obtain a good control on the error after a long time.  Nevertheless, particle simulations lead to a very good understanding of qualitative properties of solutions for aggregation equations where collisions do not happen \cite{balague_preprint,Bertozzi2015}. For the one-dimensional setting with a nonlinear dependency of the term $\nabla W \ast \rho$,  a finite volume scheme for simulating the behaviour after blow-up has been proposed in \cite{James2015NumericalMF} and  its convergence has been shown.  An energy decreasing finite volume method  for a large class of PDEs including \eqref{eq:standardmodelmacroscopic} has been proposed in \cite{carrillo2014} and a convergence result for a finite volume scheme with general measures as initial data has been shown in \cite{Carrillo2016304}. In particular, this numerical scheme leads to numerical simulations of solutions in dimension greater than one.
	
	The isotropic aggregation equation \eqref{eq:standardmodelmacroscopic} may also be modified to include linear or nonlinear diffusion terms \cite{carrillobookchapter}. While a linear diffusion term can be used to describe noise at the level of interacting particles, a nonlinear diffusion term can be used to model a system of interacting particles at the continuum level, and can be expressed by a repulsive potential. To see the latter, we consider the potential $W_\delta=W+\delta \delta_0$ for a parameter $\delta>0$ and the Dirac delta $\delta_0$, inducing an additional strongly localised repulsion. This corresponds to a PDE with nonlinear diffusion which is given by 
	\begin{align*}
	\rho_t+\nabla\cdot (\rho (-\nabla  W\ast \rho))=\delta \nabla \cdot (\rho \nabla \rho).
	\end{align*}
	More generally, adding nonlinear diffusion in \eqref{eq:standardmodelmacroscopic} results in the  class of aggregation equations
	\begin{align}\label{eq:standardmodelmacroscopicdiffusion}
	\rho_t+\nabla\cdot (\rho (-\nabla  W\ast \rho))=\delta \nabla \cdot (\rho \nabla \rho^{m-1})
	\end{align}
	with diffusion coefficient $\delta>0$ and a real exponent $m>1$. Equation \eqref{eq:standardmodelmacroscopicdiffusion} is the isotropic counterpart of \eqref{eq:macroscopiceqnonlin} for $m=2$. Of central importance for studies of \eqref{eq:standardmodelmacroscopicdiffusion} is its gradient flow formulation \cite{gradientflows} with respect to the energy
	\begin{align}\label{eq:energyisotropicdiffusion}
	\E_\delta(\rho)=\frac{1}{2}\int_{\R^d}  \rho ( W\ast \rho + \delta \rho^{m-1})\di u.
	\end{align}
	In particular, stationary states of \eqref{eq:standardmodelmacroscopicdiffusion} are critical points of the energy \eqref{eq:energyisotropicdiffusion}. The existence of global minimisers of \eqref{eq:energyisotropicdiffusion} has recently been studied in \cite{BEDROSSIAN} using techniques from the calculus of variations. While  radially symmetric and non-increasing global minimisers exist for $m>2$, the case $m = 2$ is critical and yields a global minimiser only for small enough diffusion coefficients $\delta>0$. Burger et al.\   \cite{burger2013} have shown that the threshold for $\delta$ is $\|W\|_{L^1}$ for $m=2$. Energy considerations have also been employed in \cite{Burger2008} to study the large time behaviour of solutions to \eqref{eq:standardmodelmacroscopicdiffusion} in one dimension. The existence of finite-size, compactly supported stationary states for the general power exponent $m > 1$ is investigated in \cite{Burger2014}. The uniqueness/non-uniqueness criteria are determined by the parameter $m$, with the critical power being $m=2$ \cite{Yaouniqueness}. In particular, the steady state is unique for a fixed mass for any attractive potential and $m\geq 2$. 
	
	\subsection{Contributions}
	In this work, we consider the macroscopic equations
	\eqref{eq:macroscopiceq} and \eqref{eq:macroscopiceqnonlin}. No gradient flow formulation exists in this case and   stationary solutions of the anisotropic aggregation equation generally cannot be regarded as  minimizers of an energy functional. 
	
	As a first aim of this paper, we derive equilibrium conditions for  stationary solutions  of \eqref{eq:macroscopiceq} and \eqref{eq:macroscopiceqnonlin}. Under the assumption that the stationary solutions are given by vertical line patterns, we show that the stationary solutions of the two-dimensional problem satisfy a one-dimensional equilibrium condition. This dimension reduction allows us to derive a scalar potential. We define an energy functional which depends on the scalar potential and  the diffusion coefficient $\delta$. The dimension reduction allows us to use  existing results on the stationary solutions of the anisotropic equations \eqref{eq:macroscopiceq}: minimizers of the energy functional exist and the stationary solutions of \eqref{eq:macroscopiceqnonlin} are minimisers of an  energy functional. 
	The dependence of the energy on $\delta$ can be regarded as a regularisation of the energy functional and gives rise to a sequence of  energy functionals indexed by $\delta >0$. For the sequence of energy functionals, we establish $\Gamma$-convergence  for vanishing diffusion and prove  convergence of minimisers of the regularised energy functional to minimisers of the non-regularised energy functional. 	
	
	The second aim of this paper is to investigate the dependence of the diffusion coefficient $\delta$ on stationary solutions numerically by considering an appropriate numerical scheme for the anisotropic interaction equation \eqref{eq:macroscopiceqnonlin} without gradient flow structure. The numerical scheme and its analysis is based on \cite{carrillo2014,Carrillo2016304}. The additional diffusion in the mean-field model  \eqref{eq:macroscopiceqnonlin} results in beautiful pattens which are better than the ones obtained with the particle model since too low particle numbers may result in dotted line patterns.
	
	This paper is organised as follows. In Section \ref{sec:stationarysol}, we consider stationary solutions for general underlying tensor fields first, before restricting ourselves to spatially homogeneous tensor fields  whose support is given by line patterns. For this case, we derive equilibrium conditions which can be reformulated as the minimisers of an energy functional. We show the existence of energy minimisers, and prove $\Gamma$-convergence of the regularised energies and the convergence of minimisers of the regularised energies to  minimisers of the non-regularised energy functional as the diffusion coefficient goes to zero. We consider a numerical scheme for the anisotropic, nonlocal aggregation equation with nonlinear diffusion \eqref{eq:macroscopiceqnonlin} and prove its weak convergence as the diffusion coefficient goes to zero in Section \ref{sec:numericsequilibria}. Finally, we show numerical results  in Section~\ref{sec:numresults}.

	\section{Stationary solutions}\label{sec:stationarysol}
	
	In this section, we study stationary solutions of the nonlocal aggregation equation with nonlinear diffusion \eqref{eq:macroscopiceqnonlin}. Since most applications require measure-valued solutions, we  consider  nonnegative stationary solutions of \eqref{eq:macroscopiceqnonlin}  only. 
	The stationary solutions $\rho_\infty=\rho_\infty(x,y)$, $(x,y)\in \R^2$, of \eqref{eq:macroscopiceqnonlin} satisfy
	\begin{align*}
	\nabla\cdot \left[ \rho_\infty (F(\cdot,T(x,y)) \ast \rho_\infty-\delta \nabla \rho_\infty)\right]= 0\qquad \text{a.e.\ in }\R^2,
	\end{align*}
	implying that the argument has to be constant a.e.\ in $\R^2$. Since we are interested in stationary line patterns, the stationary solution $\rho_\infty$ should satisfy $\supp{\rho_\infty}\subsetneq \R^2$ for small diffusion coefficients $\delta>0$. Hence
	it is sufficient to require
	\begin{align}\label{eq:stationaryeqorig}
	\rho_\infty ( F(\cdot,T(x,y)) \ast \rho_\infty-\delta \nabla \rho_\infty)= 0\qquad \text{a.e.\ in }\R^2,
	\end{align}
	or equivalently
	\begin{align*}
	F(\cdot,T(x,y)) \ast \rho_\infty=\delta \nabla \rho_\infty\qquad \text{on } \supp(\rho_\infty).
	\end{align*}
	In the following, we assume that the underlying tensor field $T$ is  spatially homogeneous and we study the associated stationary solutions.

	\subsection{Notation and assumptions}
	Given a spatially homogeneous tensor field, the aim of this section is to derive a scalar force and its scalar potential in one variable which  will be used to define the associated regularised and non-regularised energy functionals.  For this, we study some properties of stationary solutions first.
	
	A stationary solution  of \eqref{eq:macroscopiceqnonlin} for any spatially homogeneous tensor field $\tilde{T}$ is a coordinate transform of a stationary solution to the mean-field equation \eqref{eq:macroscopiceqnonlin} for the tensor field $T$ with $l=(1,0)$ and $s=(0,1)$ \cite{patternformationanisotropicmodel}. This motivates to study one specific spatially homogeneous tensor field in detail. In the following, we restrict ourselves to the tensor field $T$ with $l=(1,0)$ and $s=(0,1)$.
	
	For the specific tensor field $T$, the total force $F$ in \eqref{eq:totalforcenew} reduces to 
	\begin{align*}
	F((x,y),T)=\begin{pmatrix} f_l(|(x,y)|)x\\ f_s(|(x,y)|)y\end{pmatrix}
	\end{align*}
	for $(x,y)\in   \R^2$.
	We denote the components of $F$ by $F_x, F_y$, i.e.\ $F=(F_x,F_y)\in\R^2$, and	we have
	$F_x((x,y))=f_l(|(x,y)|)x$ and $F_y((x,y))=f_s(|(x,y)|) y$ for $(x,y)\in\R$  with $|(x,y)|=\sqrt{x^2+y^2}$. Since $f_s \neq f_l$, the interaction force $F$ is anisotropic. The interaction force in models of the form \eqref{eq:macroscopiceqnonlin} usually decays very fast which motivates to assume 
	\begin{align}\label{eq:forcecoeffbd}
	f_s(|(x,y)|)=f_l(|(x,y)|)=0
	\end{align}
	for $|(x,y)|\geq 0.5$
	in the following.
	
	The stationary solutions to \eqref{eq:macroscopiceqnonlin} for the  tensor field $T$ with $l=(1,0)$ and $s=(0,1)$ form vertical line patterns where the stationary solutions  are constant along the $y$-direction. Solutions on $\R^2$ which are constant along the $y$-direction are no probability measures. By restricting the domain and considering the domain $\domain=\R\times [-0.5,0.5]$ instead of $\R^2$, we can regard stationary solutions to \eqref{eq:macroscopiceqnonlin} as probability measures which are constant along the $y$-direction. Note that this assumption on the domain $\domain$ is not restrictive and by appropriate rescaling similar results can be obtained for any  domain of the form $\R\times [a,b]$ for any $a,b \in\R$ with $a<b$.

	Considering the rather unusual domain $\Omega$ can be regarded as the starting point for reducing the problem to one spatial dimension. Due to the anisotropy of the interaction force $F$, we study  solutions to \eqref{eq:macroscopiceqnonlin} on the two-dimensional space $\Omega$  which are constant in the $y$-direction. This allows us to reduce the equilibrium conditions to a one-dimensional problem. For $(x,y)\in\Omega$, we consider stationary solutions on $\domain$ of the form
	\begin{align}\label{eq:rhovertcond}
	\rho_\infty(x,y)=\rho_\infty(x,0) \quad \text{for a.e.\ } y\in [-0.5,0.5].
	\end{align} 
	The special form \eqref{eq:rhovertcond} of the stationary solutions motivates the definition of the space $\mathcal{P}_c(\domain)$ of probability measures  which are constant in the $y$-direction. We define the space $\mathcal{P}_c(\domain)$ by
	\begin{align*}
	\mathcal{P}_c(\domain)=\left\{ \rho \in L^1_+(\domain) \colon  \int_{\domain} \rho\di(x,y)=1, \enspace \rho(x,y)=\rho(x,0)  \text{ for a.e.\ } y\in [-0.5,0.5] \right\}.
	\end{align*}	
	For a consistent definition of the convolution $F\ast \rho_\infty$, we extend $F=(F_x, F_y)$  and $\rho_\infty$, defined on $\Omega$, periodically on $\R^2$ with respect to the $y$-coordinate. For any $k\in \Z$, we set $\rho_\infty(x,y+k)=\rho_\infty(x,0)$ and $F(x,y+k)=F(x,y)$  for $(x,y)\in \Omega$  implying  $F_x(x,y+k)=F_x(x,y)$ and $F_y(x,y+k)=F_y(x,y)$. This allows us to evaluate the convolution integrals $F_x \ast \rho_\infty, F_y \ast \rho_\infty$.
	For $\rho_\infty$ satisfying \eqref{eq:rhovertcond}, we have 
	\begin{align*}
	F_y \ast \rho_\infty(x,y)&=\iint \limits_{\Omega} F_y(w,z)\rho_\infty(x-w,y-z)\di (w,z)\\&=\iint \limits_{\Omega} f_s(\sqrt{w^2+z^2}) z\rho_\infty(x-w,0)\di (w,z)=0.
	\end{align*}
	Here, we used that $F_y$ is an odd function in the $y$-coordinate with $F_y(w,z)=f_s(\sqrt{w^2+z^2}) z$ for $(w,z)\in\Omega$ and $\rho_\infty$ is constant with respect to the $y$-coordinate. 
	The convolution $F_x \ast \rho_\infty$ is of the form
	\begin{align*}
	F_x\ast \rho_\infty(x,y)&=\iint \limits_{\Omega} F_x(w,z)\rho_\infty(x-w,y-z)\di (w,z)\\&=\iint \limits_{\Omega} F_x(x-w,y-z)\rho_\infty(w,z)\di (w,z)\\
	&=\int_{\R} \rho_\infty(w,0)\int_{-0.5}^{0.5} F_x(x-w,y-z)\di z\di w.
	\end{align*}	 
	Having the convolutions $F_x\ast \rho_\infty$ and $F_y\ast \rho_\infty$ at hand, we can evaluate the equilibrium condition \eqref{eq:stationaryeqorig}.  Considering the left-hand side of \eqref{eq:stationaryeqorig} as a vector, it immediately follows that its second component vanishes.  
	Since a scalar force in one variable is required for a dimension reduction, this motivates to introduce a scalar odd function $G\colon \R\to \R$ defined by 
	\begin{align}\label{eq:defintforce}
	G(x)=\int_{-0.5}^{0.5} F_x(x,z)\di z=x \int_{-0.5}^{0.5} f_l(\sqrt{x^2+z^2})\di z,
	\end{align}
	where $G(0)=0$ and $F_x(x,z)=f_l(\sqrt{x^2+z^2})x$. Note that it is not clear what the sign of $\int_{-0.5}^{0.5} f_l(\sqrt{x^2+z^2})\di z$ is. Further note that $G$ vanishes for $x\geq 0.5$ because of \eqref{eq:forcecoeffbd}.
	Due to the periodic extension of $F_x$ along the $y$-coordinate, we have  $G(x)=\int_{-0.5}^{0.5} F_x(x,y-z)\di z$ for any $y\in [-0.5,0.5]$.
	Hence, there exists an interaction  potential $W\colon \R\to \R$  which is even and satisfies
	\begin{align}\label{eq:defG}
	G=-W'.
	\end{align}

	The one-dimensional potential $W$ and the assumption on the stationary solution $\rho_\infty$ of the form \eqref{eq:rhovertcond} imply that the problem is in fact one-dimensional. It is sufficient to consider $\rho_\infty=\rho_\infty(x)$ for $x\in \R$ and the space $\mathcal P_c$ reduces to 
	\begin{align*}
	\mathcal{P}_c(\R)=\left\{ \rho \in L^1_+(\R) \colon  \int_{\R} \rho\di x=1  \right\}.
	\end{align*} 
	Using the potential $W$, we define  the energy functional
	\begin{align}\label{eq:energy1dnonreg}
	\mathcal{E}(\rho_\infty)=\frac{1}{2}  \int_{\R} \rho_\infty ( W\ast \rho_\infty)\di x
	\end{align}
	in one spatial dimension
	where $W\ast \rho_\infty$ is the convolution in one coordinate, i.e.\
	\begin{align}\label{eq:convolutionW}
	W\ast \rho_\infty(x)=\int_\R W(x-w)\rho_\infty(w)\di w
	\end{align}
	for $x\in \R$.
	The associated  equilibrium condition is given by
	\begin{align}\label{eq:stationaryeqhomunreg}
	\rho_\infty \partial_x ( W \ast \rho_\infty)= 0\quad \text{a.e.\ in }\R.
	\end{align}
	The regularisation of the energy $\E$ is defined as
	\begin{align}\label{eq:energy1d}
	\mathcal{E}_\delta(\rho_\infty)=\frac{1}{2}  \int_{\R} \rho_\infty ( W\ast \rho_\infty+\delta \rho_\infty)\di x.
	\end{align}
	The associated equilibrium condition is given by
	\begin{align}\label{eq:stationaryeqhom}
	\rho_\infty \partial_x ( W \ast \rho_\infty+\delta\rho_\infty)= 0\qquad \text{a.e.\ in }\R.
	\end{align}
	
	One-dimensional conditions of the form \eqref{eq:stationaryeqhom} have already been studied in the literature. Stationary solutions are considered via energy minimisation in \cite{burger2013} and we state the result in Proposition \ref{prop:energymin}. The properties of stationary solutions $\rho_\infty$ have been studied in \cite{burger2013} for purely attractive potentials but the results in fact also hold in our setting. For completeness, we state these results in Section~\ref{sec:propsolution} where stationary solutions on $\R$ are considered. For the comparison of analytical and numerical results, we are interested in stationary solutions on the torus $\mathbb{T}^2$ which may also be observed in numerical simulations. We investigate  conditions for stationary solutions on $\mathbb{T}^2$ in Section \ref{sec:torus}. 
	The analytical results require rather relaxed conditions on the potential $W$:
	\begin{assumption}\label{ass:W}
		For the interaction potential $W$ satisfying \eqref{eq:defG}, we require 
		\begin{assump}\label{ass:W0}
			$W$ is even, i.e.\ $W(x)=W(-x)$.
		\end{assump} 
		\begin{assump}\label{ass:W1}
			$W$ is continuously differentiable. 
		\end{assump}
		\begin{assump}\label{ass:W3}
			$W'(x)=0$ for $|x|\geq 0.5$.
		\end{assump} 
		\begin{assump}\label{ass:W4}
			There exist $\bar{\delta}>0$ and a measure $\rho_\infty\in \mathcal{P}_c(\domain)$ such that $\E_{\bar{\delta}}(\rho_\infty)\leq 0$. 
		\end{assump} 
	\end{assumption}
	
	Assumption \ref{ass:W} is basically an assumption for $f_l$. Since the force coefficient $f_l$ is short-range repulsive, long-range attractive, we have $f_l(0)>0$, there exists $x\in(0,0.5)$ such that $f_l(x)<0$ and $f_l(x)=0$ for $x\geq 0.5$.

	\begin{remark}
		Note that assumptions \assref{ass:W0},  \assref{ass:W1} and  \assref{ass:W3} are rather relaxed conditions and allow us to consider a  general class of interaction potentials, including the one that can be derived from $G$ based on $F_x$ in \eqref{eq:defintforce}. 	In particular, the interaction potential $W$ 
		is bounded. Besides,  the energy $\E_\delta\colon \mathcal{P}_c(\R)\to \R$ in \eqref{eq:energy1d} is weakly lower semi-continuous with respect to weak convergence of measures. 
		Assumption \assref{ass:W4} is required for establishing the existence of minimisers of the energy $\E_\delta$ in \eqref{eq:energy1d}. By \assref{ass:W4}, there exists a measure $\rho_\infty\in \mathcal{P}_c(\R)$ 
		such that $\E_{\delta}(\rho_\infty)\leq 0$ for all $0\leq \delta \leq \bar{\delta}$. 
	\end{remark}

	\begin{remark}
		Assumption \assref{ass:W4} is not restrictive since  $W$ satisfying Assumption \ref{ass:W} is only given up to an additive constant by \eqref{eq:defG}. We can choose the additive constant in such a way that the boundedness of $W$ in \assref{ass:W3} guarantees \assref{ass:W4}.
		Examples for $\rho_\infty$ satisfying \assref{ass:W4} include mollified delta distributions and indicator functions of the form $\rho_\infty=\frac{1}{|Q_W|}\chi_{Q_W}$ where $Q_W=[-x_W/2,x_W/2]$ for some $x_W>0$. In particular, $\rho_\infty$ may have compact, connected support.
		We recall that Assumption \assref{ass:W4} is equivalent to the existence of minimizers of $\E_\delta$. In addition, it is  equivalent to $\E_\delta$ not being $H$-stable \cite{Canizo2015,Simione2015}. The notation of $H$-stability is important in statistical mechanics. A system of interacting particles has a macroscopic thermodynamic behaviour provided mass is not accumulated on bounded regions as the number of particles goes to infinity. Such potentials  are called $H$-stable. We say a  potential  $W$ is $H$-stable if  there exists $B\in \R$ such that for all $N$ and for all sets of $N$ distinct points $\{x_1,\ldots,x_N\}$ in $\R^2$ it holds 
		$$\frac{1}{N^2}\sum_{1\leq i<j\leq N} W(x_i-x_j) \geq -\frac{1}{N}B.$$
		The $H$-stability of a potential is equivalent to $\E(\rho)\geq 0$ for any probability measure $\rho$.
	\end{remark}
	
	\subsection{Equilibrium conditions}
	In this section, we consider the equilibrium condition \eqref{eq:stationaryeqhom}.
	Since we are only interested in minimisers of the interaction energy $\E_\delta$, we require
	\begin{align}\label{eq:stationaryeq}
	W \ast \rho_\infty+\delta \rho_\infty=C \qquad  \text{in each connected component of } \supp(\rho_\infty)
	\end{align}
	for some constant $C\in\R$. 	
	By multiplying \eqref{eq:stationaryeq} by $\rho_\infty$ and integrating over $\supp(\rho_\infty)$, we obtain
	\begin{align*}
	\int_{\supp(\rho_\infty)} \rho_\infty (W \ast \rho_\infty) \di x+\delta \int_{\supp(\rho_\infty)} \rho_\infty^2\di x=C,
	\end{align*}
	where  the unit mass of $\rho_\infty$ was used. In particular, this shows that $C=C(\delta)\in \R$ is uniquely determined and the integral equation \eqref{eq:stationaryeq}  may be expressed in the equivalent fixed point form 
	\begin{align}\label{eq:rhosteadystate}
	\rho_\infty=\frac{(C-W\ast \rho_\infty)_+}{\int_{\supp(\rho_\infty)} (C-W\ast \rho_\infty)_+ \di x}.
	\end{align}
	Clearly, the fixed point form is consistent with \eqref{eq:rhovertcond} and the dependence of $\rho_\infty$ on $\delta$ follows from $C=C(\delta)$.
	
	Non-trivial stationary states with purely attractive potentials in the set $L^2(\R^d)\cap \mathcal{P}(\R^d)$ with $d\geq 1$ are considered in \cite{burger2013}. The authors show that minimisers of the energy functional \eqref{eq:energy1d} are sufficient for solving the equilibrium conditions.  Analogously, one can  show the same result in our setting for potentials satisfying Assumption \ref{ass:W} and stationary states in the space $L^2(\R)\cap \mathcal{P}_c(\R)$. 
	In particular, a minimiser of the energy functional \eqref{eq:energy1d} is sufficient for solving \eqref{eq:stationaryeqhom}.
	\begin{proposition}[Stationary solutions via energy minimisation] \label{prop:energymin}
		Let $\rho_\infty\in L^2(\R)$ be a minimiser of the energy functional \eqref{eq:energy1d} on $\mathcal{P}_c(\R)$.
		Then, $\rho_\infty$ satisfies \eqref{eq:stationaryeqorig}.
	\end{proposition}

	\subsection{Existence and convergence of minimisers}\label{sec:minimizers}
	Motivated by Proposition \ref{prop:energymin}, we consider the energy functionals $\E$ and $\E_\delta$, defined in \eqref{eq:energy1dnonreg} and \eqref{eq:energy1d}.
	For the existence  and convergence of minimisers, we have to verify that an energy minimising sequence is precompact in the sense of weak convergence of measures, and prove a $\Gamma$-convergence result. For this, we use Lions' concentration compactness lemma for probability measures \cite{lions1984}, \cite[Section 4.3]{Struwe2000} and reformulate it to our setting. 
	
	\begin{lemma}[Concentration-compactness lemma for measures]\label{lem:lions}
		Let $\{\rho_n\}_{n\in\N}\subset \mathcal{P}_c(\R)$. Then, there exists a subsequence $\{\rho_{n_k}\}_{k\in\N}$ satisfying one of the three following possibilities:
		\begin{enumerate}
			\item (tightness up to transition) There exists $z_k\in\R$ such that for all $\epsilon>0$ there exists $R>0$ satisfying
			\begin{align*}
			\int_{B_R(z_k)\cap \R} \di \rho_{n_k} (x)\geq 1-\epsilon\quad \text{for all }k;
			\end{align*}
			\item\label{lem:lions2} (vanishing) 
			\begin{align*}
			\lim_{k\to \infty} \sup_{z\in\R} \int_{B_R(z)\cap \R} \di \rho_{n_k} (x)=0\quad \text{for all }R>0;
			\end{align*}
			\item (dichotomy) There exists $\alpha\in(0,1)$ such that for all $\epsilon>0$ there exists $R>0$ and a sequence $\{z_k\}_{k\in\N}\subset \R$ with the following property: 
			
			Given any $R'>R$ there are nonnegative measures $\rho_k^1$ and $\rho_k^2$ such that
			\begin{align*}
			&0\leq \rho_k^1+\rho_k^2\leq \rho_{n_k},\\
			&\supp(\rho_k^1)\subset B_R(z_k)\cap \R,\\
			&\supp(\rho_k^2)\subset \R \backslash B_{R'}(z_k),\\
			&\limsup_{k\to \infty}\bl \left| \alpha -\int_{\R} \di \rho_k^1(x)\right|+\left| (1-\alpha) -\int_{\R} \di \rho_k^2(x)\right|\br \leq \epsilon.
			\end{align*}
		\end{enumerate}	
	\end{lemma}
	
	For proving the existence of minimisers of the energy functional \eqref{eq:energy1d}, one can use the direct method of the calculus of variations and Lemma \ref{lem:lions} to eliminate the cases `vanishing' and `dichotomy' of an energy minimising sequence. The proof of the existence of minimisers of the regularised energy $\E_\delta$ in \eqref{eq:energy1d} is very similar to the one for the non-regularised energy $\E$,
	provided in \cite[Theorem 3.2]{Simione2015}:
	
	\begin{proposition}[Existence of minimisers]
		Suppose $W$ satisfies assumptions \assref{ass:W0},  \assref{ass:W1} and \assref{ass:W3}. Then, the regularised energy $\E_\delta$ in \eqref{eq:energy1d} has a global minimiser in $\mathcal{P}_c({\R})$ if and only if it satisfies \assref{ass:W4}.  The non-regularised energy $\E$ in \eqref{eq:energy1dnonreg} has a global minimiser in $\mathcal{P}_c({\R})$ if and only if  \assref{ass:W4} is satisfied for $\E$. 
	\end{proposition}

	\begin{theorem}[$\Gamma$-convergence of regularised energies]\label{th:gammaconv}
		Suppose that $W$ satisfies \assref{ass:W0}, \assref{ass:W1} and \assref{ass:W3}. The sequence of regularised energies $\{\E_\delta\}_{\delta>0}$ $\Gamma$-converges to the energy $\E$ with respect to the weak convergence of measures. That is, 
		\begin{itemize}
			\item (Liminf) For any $\{ \rho_\delta \}_{\delta>0}\subset \mathcal{P}_c({\R})$ and $\rho\in \mathcal{P}_c({\R})$ such that $\rho_\delta$  converges weakly to $\rho$ as $\delta \to 0$, we have
			\begin{align*}
			\liminf_{\delta \to 0} \E_\delta (\rho_\delta)\geq \E(\rho).
			\end{align*}	
			\item (Limsup) For any $\rho \in \mathcal{P}_c({\R})$ there exists a sequence $\{\rho_\delta\}_{\delta>0}\in\mathcal{P}_c({\R})$ such that $\rho_\delta$  converges weakly to $\rho$ as $\delta\to 0$ and
			\begin{align*}
			\limsup_{\delta \to 0} \E_\delta (\rho_\delta)\leq \E(\rho).
			\end{align*}	
		\end{itemize}	
	\end{theorem} 
	
	\begin{proof}
		Step 1 (Liminf): Since $W$ is lower semi-continuous and bounded from below, the weak lower semi-continuity of the first term in the energy functional $\E_\delta$ in \eqref{eq:energy1d}  follows from the Portmanteau Theorem \cite[Theorem 1.3.4]{vaartwellner96book}, i.e.
		\begin{align*}
		\liminf_{\delta\to 0} \frac{1}{2} {\int_{\R} \rho_\delta ( W\ast \rho_\delta)\di x} \geq \frac{1}{2}  {\int_{\R} \rho ( W\ast \rho)\di x}.
		\end{align*}
		Together with
		\begin{align*}
		\liminf_{\delta \to 0}\frac{\delta}{2} { \int_{\R} \rho_\delta^2\di x}\geq 0,
		\end{align*}
		the liminf inequality immediately follows.
		
		Step 2 (Limsup): Let $\mu\in\mathcal{P}_c({\R})$ be given, let $$\phi(x)=\frac{1}{\sqrt{4\pi}}\exp\left(-\frac{|x|^2}{4}\right)$$ denote the one-dimensional heat kernel and define $$\phi_\delta(x)=\frac{1}{\sqrt{\delta}}\phi\left(\frac{x}{\sqrt{\delta}}\right).$$
		Note that $\phi\in C^\infty({\R})$, ${\phi(x)=\phi(-x)}$ for all ${x\in\R}$, and $${\int_\R\phi \di x=1}.$$ 
		In particular, $|\phi_\delta|\leq \frac{C_\phi}{\sqrt{\delta}}$ where $C_\phi$ denotes the bound of $\phi$.
		We define the measure $\rho_\delta:=\phi_\delta \ast \rho$ which converges weakly to $\rho$ in $\mathcal{P}_c({\R})$.  Note that 
		\begin{align*}
		\delta \int_\R \rho_\delta^2 \di x\leq C_\phi\sqrt{\delta} \int_\R\rho_\delta \di x=C_\phi\sqrt{\delta}\to 0\qquad \text{as }\delta \to 0.
		\end{align*}
		Due to the continuity of $W$,  the term ${- \int_{\R} \rho( W\ast \rho)\di x}$ is weakly lower semi-continuous and
		\begin{align*}
		{\limsup_{\delta\to 0} \frac{1}{2}  \int_{\R} \rho_\delta ( W\ast \rho_\delta)\di x \leq \frac{1}{2}  \int_{\R} \rho ( W\ast \rho)\di x,}
		\end{align*}
		resulting in the limsup inequality.
	\end{proof}
	
	\begin{theorem}[Convergence of minimisers]
		Suppose that $W$ satisfies \assref{ass:W0}, \assref{ass:W1} and \assref{ass:W3}.	For any $\bar{\delta}>0$ sufficiently small, suppose that $\E_{\bar{\delta}}$ satisfies \assref{ass:W4}  and let $\rho_\delta \in \mathcal{P}_c({\R})$ be a minimiser of the energy $\E_\delta$ in \eqref{eq:energy1d} for all $0<\delta\leq \bar{\delta}$. Then, there exists $\rho\in\mathcal{P}_c({\R})$ such that, up to a subsequence and translations, $\rho_\delta$ converges weakly to $\rho$ as $\delta \to 0$, and $\rho$ minimises the energy $\E$ over $\mathcal{P}_c({\R})$.
	\end{theorem} 	  
	
	\begin{proof} 		
		Let $\{\rho_\delta\}_{\delta>0}\subset \mathcal{P}_c({\R})$ be a sequence of minimisers of $\E_\delta$. For $\bar{\delta}>0$ sufficiently small, we may assume that $\E_\delta(\rho_\delta)\leq 0$ for all $0<\delta\leq \bar{\delta}$ since $\rho_\delta$ minimises $\E_\delta$.  As in \cite[Theorem 3.2]{Simione2015} one can eliminate the cases `vanishing' and `dichotomy' in Lemma \ref{lem:lions}, implying that there exists a subsequence $\{\rho_{\delta_k}\}_{k\in\N}$ satisfying `tightness up to translation', i.e.\ there exists $z_k\in{\R}$ such that for all $\epsilon>0$ there exists $R>0$ satisfying
		\begin{align*}
		\int_{B_R(z_k)\cap \R} \di \rho_{\delta_k}{(x)}\geq 1-\epsilon\quad \text{for all }k.
		\end{align*}
		We define $\tilde{\rho}_{\delta_k}:=\rho_{\delta_k}(\cdot-z_k)$ and hence $\{\tilde{\rho}_{\delta_k}\}_{k\in\N}$ is tight. Since $\E_{\delta_k}(\rho_{\delta_k})=\E_{\delta_k}[\tilde{\rho}_{\delta_k}]$,  $\{\tilde{\rho}_{\delta_k}\}_{k\in\N}$ is also a sequence of minimisers of $\E_{\delta_k}$ and by Prokhorov's Theorem (cf.\ \cite[Theorem 4.1]{Billingsley}) there exists a further subsequence $\{\tilde{\rho}_{\delta_k}\}_{k\in\N}$, not relabelled, such that $\tilde{\rho}_{\delta_k}$ converges weakly to some measure $\rho\in\mathcal{P}_c(\R)$ as $k\to\infty$.
		
		For showing that the measure $\rho$ minimises the energy functional $\E$, we consider an arbitrary measure $\mu\in\mathcal{P}_c(\R)$. By the limsup inequality in Theorem \ref{th:gammaconv}, there exists a sequence $\{\mu_{\delta_k}\}_{k\in\N}$ which  converges weakly to $\mu$ as $k\to\infty$ such that 
		\begin{align*}
		\limsup_{k   \to \infty} \E_{\delta_k} (\mu_{\delta_k})\leq \E(\mu).
		\end{align*}	
		Together with the liminf inequality in Theorem \ref{th:gammaconv}, this yields
		\begin{align*}
		\lim_{k \to \infty} \E_{\delta_k} (\mu_{\delta_k})= \E(\mu).
		\end{align*}
		Since the sequence of measures $\tilde{\rho}_{\delta_k}$ is a minimising sequence of $\E_{\delta_k}$ which  converges weakly to $\rho$, we obtain, again by the liminf inequality,
		\begin{align*}
		\E(\rho)\leq \liminf_{k \to \infty} \E_{\delta_k} (\tilde{\rho}_{\delta_k})\leq \liminf_{k \to \infty} \E_{\delta_k} (\mu_{\delta_k})= \E(\mu).
		\end{align*}
	\end{proof}  
	
	Note that each local minimiser $\rho_\infty$  of $\E_\delta$ is a steady state and satisfies the equilibrium condition \eqref{eq:stationaryeqhom}. To see this, note that the Euler-Lagrange conditions for minimisers \cite[Proposition 2.4]{Carrillo2019} state that  for each connected component $A_i$ of $\supp(\rho_\infty)$ there exists $C_i\in \R$ such that
	\begin{align*}
	W \ast \rho_\infty+\delta \rho_\infty&=C_i \qquad  \text{a.e.\ on } A_i,\\
	W \ast \rho_\infty+\delta \rho_\infty&\geq C_i \qquad  \text{a.e.\ on } \R.
	\end{align*}
	Since $\supp(\rho_\infty)$ is connected in our setting, see Theorem \ref{th:existencesteadystate} below, this implies that \eqref{eq:stationaryeq} is fulfilled, implying that $\rho_\infty$ is of the form \eqref{eq:rhosteadystate}. In particular, $\partial_x \rho$ is well-defined and condition \eqref{eq:stationaryeqhom} holds.

	\subsection{Properties of stationary solutions}\label{sec:propsolution}
	
	The sign of the odd function $G$, defined by $G(x)=\int_{-0.5}^{0.5} F_x(x,z)\di z$ in \eqref{eq:defintforce}, is not clear for the force $F_x$ in the K\"ucken-Champod model, see \cite{During2017} or Section \ref{sec:intro} for the precise definition of the force coefficients. Since  $G=-W'$, $W$ is only determined up to an additive constant. Due to the boundedness of $W$ in Assumption \ref{ass:W}, we can chose the additive constant such that
	\begin{align}\label{eq:propertiesWnew}
	W(x)\leq 0 \enspace \text{for all}\enspace |x|\geq 0.
	\end{align}
	In particular, the assumptions on the potential $W$ for the one-dimensional results in \cite{burger2013} are satisfied for \eqref{eq:propertiesWnew} and the results also hold for our setting. For completeness, we state  results on properties of stationary solutions, proven in \cite{burger2013}, below: 
	
	\begin{corollary}
		Let $\delta>0$ be given. 
		\begin{itemize}
			\item If $\delta \geq \|W\|_{L^1}$, there exists no stationary solution  $\rho_\infty$ in $L^2\cap \mathcal{P}_c(\R)$ satisfying \eqref{eq:stationaryeqhom}.
			\item If  $\delta < \|W\|_{L^1}$,  there exists a minimiser $\rho_\infty \in L^2\cap \mathcal{P}_c(\R)$ of the energy functional \eqref{eq:energy1d}  which is  symmetric in $x$,  non-increasing on $x\geq 0$, and satisfies \eqref{eq:stationaryeqhom}.
		\end{itemize}	
	\end{corollary}	
	
	To relate the cases $\delta<\|W\|_{L^1}$ and $\delta\geq\|W\|_{L^1}$ to  \assref{ass:W4} note that
	\begin{align*}
	- \int_{\R} \rho_\infty (W\ast \rho_\infty)\di x\leq \|W\|_{L^1} \int_{\R} \rho_\infty^2 \di x
	\end{align*}
	by Young's convolution inequality and property \eqref{eq:propertiesWnew} of $W$, implying
	\begin{align*}
	\E_\delta(\rho_\infty)=\frac{1}{2}  \int_{\R} \rho_\infty ( W\ast \rho_\infty+\delta \rho_\infty)\di x\geq \frac{\delta- \|W\|_{L^1} }{2}  \int_{\R} \rho_\infty^2 \di x.
	\end{align*}
	A necessary condition for \assref{ass:W4} is given by $\delta \leq \|W\|_{L^1} $. 
	
	\begin{proposition}\label{prop:relationepsl}
		For any given $L > 0$ there exists a unique symmetric function $\rho_\delta \in C^2([-L,L])$ with unit mass 
		and $\frac{\di }{\di x}\rho_\delta(x) \leq 0$ for $x \geq 0$
		such that $\rho_\delta$ solves \eqref{eq:stationaryeq} for some $\delta = \delta(L) > 0$ where $C=C(\delta)$ in \eqref{eq:stationaryeq} satisfies $C=2\E_{\delta}(\rho_\delta)$. Such a function $\rho_\delta$ also satisfies $\frac{\di^2 }{\di x^2} \rho_\delta(0) < 0$.
		Moreover, $\delta(L)$ is the largest eigenvalue of the compact operator 
		\begin{align*}
		\mathcal{W}_L[\rho_\delta](x):&=\int_0^L \rho_\delta(w)  \bigg( W(x-w)+W(x+w)-W(L-w)-W(L+w)\bigg) \di w
		\end{align*}
		on the Banach space  $$\mathcal{Y}_L:=\{\rho_\delta\in C([0,L]\times [-0.5,0.5]) \colon \rho_\delta(L,y)=0 \text{ for all }y\in[-0.5,0.5] \}.$$ 
		The simple eigenvalue $\delta(L)$  is uniquely determined as a function of $L$ with the following properties:
		\begin{enumerate}
			\item $\delta(L)$ is continuous and strictly increasing with respect to $L$, 
			\item $\lim_{L\to +\delta} \delta(L) = \|W\|_{L^1}$,
			\item $\delta(0) = 0$.
		\end{enumerate}
	\end{proposition}
	
	\begin{theorem}\label{th:existencesteadystate}
		Let $0<\delta < \|W\|_{L^1}$. Then, there exists a unique $\rho_\delta \in  L^2\cap \mathcal{P}_c(\R)$ with unit mass and zero centre of mass  such that \eqref{eq:stationaryeqhom} is  satisfied.
		Moreover,  
		\begin{itemize}
			\item $\rho_\delta$ is symmetric in $x$ and monotonically decreasing on $x > 0$,
			\item $\rho_\delta\in C^2(\supp (\rho_\delta))$,
			\item  $\supp(\rho_\delta)$ is a bounded, connected set in $\R$,
			\item $\rho_\delta$ has a global maximum at $x=0$, and $\frac{\di ^2}{\di x^2}\rho_\delta(0)<0$,
			\item $\rho_\delta$ is the global minimiser of the energy $\E_\delta$ in \eqref{eq:energy1d}.
		\end{itemize}
	\end{theorem}	
	
	\subsection{Stationary solutions on the torus}\label{sec:torus}
	We studied the steady states on $\Omega=\R\times [-0.5,0.5]$ in the previous subsections, equivalent to the one-dimensional problem on $\R$. In this subsection, we investigate the steady states on the two-dimensional unit torus $\mathbb{T}^2$, or equivalently, the unit square $[-0.5,0.5]^2$ with periodic boundary conditions. Here, we focus on steady states which exist under perturbation of the potential. This is also helpful for comparing the  analytical results to the numerical simulations. 
	
	As on the domain $\Omega=\R \times [-0.5,0.5]$, we require that  minimisers  $\rho_\infty$ of the energy functional $\E_\delta$ in \eqref{eq:energy1d} are constant in $y$-direction and of the form $\rho_\infty(x,y)=\rho_\infty(x,0)$ for all $y\in[-0.5,0.5]$ with zero centre of mass. This assumption allows us to study the associated one-dimensional problem.
	
	While we have seen in Theorem \ref{th:existencesteadystate} that steady states on $\Omega=\R\times [-0.5,0.5]$ have a connected support, we show in this subsection that steady states on the unit torus may not have connected support and may be composed of finitely many stripes of equal width and equal distances between each other. To see this, let us consider minimisers of the non-regularised energy $\E$ in \eqref{eq:energy1dnonreg}. We assume that the associated steady state $\rho_\infty$ is of the form
	\begin{align}\label{eq:sumdelta}
	\rho_\infty(x)=\frac{1}{n} \sum_{k=1}^n\delta_{x_k}(x)
	\end{align}
	for $x_1,\ldots,x_n\in(-0.5,0.5)$ with $x_1<\ldots<x_n$. 
	The equilibrium condition \eqref{eq:stationaryeqhomunreg} implies $W'\ast \rho_\infty=0$ a.e.\ on $\supp(\rho_\infty)$, i.e.		
	\begin{align}\label{eq:multiplelinescond}
	\sum_{j=1}^n W'(x_k-x_j)=\sum_{\substack{j=1 \\ j\neq k}}^n W'(x_k-x_j) =0
	\end{align}
	for all $k\in\{1,\ldots,n\}$ where we used that $W'(0)=-G(0)=0$.
	
	Let us suppose first that we have an odd number $n$ of stripes.  
	Condition \eqref{eq:multiplelinescond} is  satisfied for any potential $W$ in Assumption \ref{ass:W} for equidistant points $x_1,\ldots,x_n$ in \eqref{eq:sumdelta} with 
	\begin{align}\label{eq:diraclocodd}
	x_k=\frac{k}{n}-\frac{n+1}{2n}, \qquad k=1,\ldots,n,
	\end{align} 
	since $W'(x)=-W'(-x)$ for $x\in\R$. In this case, all conditions in \eqref{eq:multiplelinescond} are equivalent to each other. If  \eqref{eq:diraclocodd} is not fulfilled, one can construct potentials in Assumption  \ref{ass:W} such that  \eqref{eq:multiplelinescond} is satisfied for all $k=1,\ldots,n$, but it may be cumbersome to satisfy the $n$ conditions in \eqref{eq:diraclocodd}. In general, there exist perturbations of these constructed potentials so that the perturbed potential does no longer satisfy the steady state condition \eqref{eq:multiplelinescond}, suggesting that non-equidistant line patterns are unstable steady states in general. Since any potential in Assumption \ref{ass:W} satisfies \eqref{eq:multiplelinescond} for equidistant stripes at positions \eqref{eq:diraclocodd}, any perturbation of the potential leads to the same steady state of equidistant lines. This suggests that it may be a stable steady state. The steady state $\rho_\infty$ is of the form \eqref{eq:sumdelta} with zero centre of mass satisfying \eqref{eq:multiplelinescond} and consisting of an odd number  $n$ of parallel equidistant lines   at locations $x_k$ in \eqref{eq:diraclocodd}. In particular, the single straight vertical line with zero centre of mass is included in the property of locations $x_k$ in \eqref{eq:diraclocodd}.

	For an even number $n$ of lines, we can proceed in a similar way as above. Condition \eqref{eq:multiplelinescond} implies that for stable steady states consisting of an even number of lines the property  $W'(-0.5)=W'(0.5)=0$ is required  in addition to equidistant lines at locations $x_k$ in \eqref{eq:diraclocodd}. Note that $W'(-0.5)=W'(0.5)=0$ is equivalent to $f_l(0.5)=0$  for the force coefficient $f_l$ in the definition of the force $F_x((x,y))=f_l(|(x,y)|)(x,y)$ for $(x,y)\in \R^2$, satisfied by \eqref{eq:forcecoeffbd}.
	
	In the following, we  study more general stable steady states $\rho_\infty$ in the form of line patterns where we suppose  that $\supp(\rho_\infty)$ consists of $n$ connected components $M_k, k=1,\ldots,n.$  Motivated by the results in Theorem \ref{th:existencesteadystate}, we consider the regularisation parameter $\delta>0$ and suppose  that the support  of $\rho_\infty$ has a positive measure. 	
	The regularised steady state condition \eqref{eq:stationaryeqhom} is given by
	\begin{align*}
	W' \ast \rho_\infty+\delta \partial_x\rho_\infty= 0\qquad \text{a.e.\ in }\supp(\rho_\infty).
	\end{align*}
	We assume that the connected components $M_k$ are intervals which are ordered in such a way that $M_k=[u_k,v_k]$ for some $u_k<v_k$ for $k=1,\ldots,n$ and $v_k<u_{k+1}$ for $k=1,\ldots,n-1$. 
	Further, we assume that each connected component $M_k$ is of the same size, i.e.\ $v_k-u_k=v_{k+1}-u_{k+1}$ for $k=1,\ldots,n-1$.
	We also assume that $\rho_\infty$ is symmetric in $x$ on $M_k$ for $k=1,\ldots,n$ as in Theorem \ref{th:existencesteadystate}. Note that for measures with zero centre of mass, we can assume without loss of generality that $\rho_\infty(-0.5)=\rho_\infty(0.5)=0$. Due to the symmetry of $\rho_\infty$ in $x$ on each $M_k$, one can see as in the unregularised case that equidistant lines as in \eqref{eq:diraclocodd} are required for stable steady states independent of the choice of the potential $W$ in Assumption \ref{ass:W}. This suggests to consider $n$ connected components $M_k$ which form
	$n$ equidistant lines, given by
	\begin{align}\label{eq:translationsupport}
	M_k=M_j+\frac{k-j}{n}, \qquad j,k\in\{1,\ldots,n\}. 
	\end{align}
	In particular, the associated steady state $\rho_\infty$ is periodic in $x$ with period $\frac{1}{n}$. As before,  $W'(-0.5)=W'(0.5)=0$ has to be satisfied for $n$ even which clearly holds due to the equivalence  to $f_l(0.5)=0$ in \eqref{eq:forcecoeffbd}. 
	In particular, this shows that the energy functionals $\E_\delta$ and $\E$ for probability measures defined on the torus $\mathbb{T}^2$ may have multiple local minimisers due to the dependence on $n$. 
	The support of these minimisers may not be connected and may consist of a finite number of connected components of equal size, satisfying the translation property \eqref{eq:translationsupport}. Besides, symmetry in $x$ on each connected component $M_k$ is required for minimisers, implying the periodicity of minimisers in $x$.

	\section{Numerical scheme and its convergence}\label{sec:numericsequilibria}	
	
	\subsection{Numerical methods}	
	For the numerical simulations, we consider the positi-
	vity-preserving finite-volume method for nonlinear equations with gradient structure proposed in \cite{carrillo2014} for isotropic interaction equations \eqref{eq:standardmodelmacroscopic}.  We consider the domain $\R^2$ and extend the scheme \cite{carrillo2014}  to the anisotropic interaction equations with or without diffusion in \eqref{eq:macroscopiceqnonlin} or \eqref{eq:macroscopiceq}, respectively. 
	This is achieved by replacing $-\nabla W$ by $F(\cdot,T)$, requiring additional care in calculating the term $(F(\cdot,T(x,y)) \ast \rho(t,\cdot))(x,y)$ for $(x,y)\in\R^2$ efficiently.
	
	In two spatial dimensions, we consider a Cartesian grid, given by $x_i=i\Delta x$ and $y_j=j\Delta y$ for $i,j\in\Z$. Let $C_{ij}$ denote the cell of the spatial discretisation $C_{ij}=[x_i,x_{i+1})\times [y_j,y_{j+1})$, and let the time discretisation be given by 
	$t_n=\sum_{i=0}^{n-1}\Delta t^i$  for $n\in\N$ with time steps $\Delta t^i$. Let $\rho_{ij}^n$ denote the approximation of the solution $\rho(t_n,x_i,y_j)$ to the anisotropic nonlocal interaction equation with diffusion \eqref{eq:macroscopiceqnonlin} with  initial condition $\rho|_{t=0}  = \rho^{in}$ in $\R^2$ for a given probability measure $\rho^{in}$. Note that \eqref{eq:macroscopiceqnonlin} can be rewritten as 
	\begin{align*}
	\partial_t \rho +\nabla \cdot \left(\rho u_\rho \right) =\delta  \nabla\cdot (\rho \nabla \rho)
	\end{align*}
	where $u_\rho$ is defined in \eqref{eq:velocityfield} with 
	\begin{align*}
	|u_\rho(t,x,y)| \leq f
	\end{align*} 
	for the uniform bound $f$ of $F$. Assuming that $\rho^{in}\in\mathcal{P}_2(\R^2)$ where $\mathcal{P}_2(\R^2)$ denotes the space of probability measures with finite second order moment, we define its discretisation
	\begin{align}\label{eq:initialconddiscrete}
	\rho_{ij}^0=\frac{1}{\Delta x \Delta y}\iint \limits_{C_{ij}} \rho^{in} \di (x, y) \geq 0
	\end{align}
	for $(i,j)\in\Z^2$. Since $\rho^{in}$ is a probability measure, the total mass of the system is $\sum_{i,j} \rho^0_{ij}\Delta x\Delta y=1$ initially. Given an approximating sequence $\{\rho_{ij}^n\}_{i,j}$ at time $t_n$, we consider the scheme
	\begin{align}\label{eq:numericalscheme}
	\begin{split}
	\rho_{ij}^{n+1}&=\rho_{ij}^n-\frac{\Delta t^n}{\Delta x} \bl {(u_x)}^n_{i+1/2,j} \rho^n_{i+1/2,j} - {(u_x)}_{i-1/2,j}^n \rho^n_{i-1/2,j}\br\\&\quad-\frac{\Delta t^n}{\Delta y} \bl {(u_y)}^n_{{i,j+1/2}} \rho^n_{i,j+1/2} - {(u_y)}^n_{{i,j-1/2}} \rho^n_{i,j-1/2}\br\\&\quad+\frac{\Delta t^n}{2\Delta x}f\bl \rho_{i+1, j}^n-2\rho_{i j}^n+\rho_{i-1, j}^n\br+\frac{\Delta t^n}{2\Delta y}f\bl \rho_{i,j+1}^n-2\rho_{i j}^n+\rho_{i, j-1}^n\br\\&\quad+ \frac{\delta\Delta t^n}{2(\Delta x)^2}\bl(\rho_{i+1,j}^n)^2-2\bl\rho_{ij}^n\br^2+\bl\rho_{i-1,j}^n\br^2\br\\&\quad+\frac{\delta \Delta t^n}{2(\Delta y)^2}\bl(\rho_{i,j+1}^n)^2-2\bl\rho_{ij}^n\br^2+\bl\rho_{i,j-1}^n\br^2\br
	\end{split}
	\end{align}
	for the uniform bound $f$ of the force $F$ and parameter $\delta>0$. Here, we use the notation
	\begin{align*}
	\rho_{i+1/2, j}&=\frac{\rho_{ij}+\rho_{i+1,j}}{2},\qquad  \rho_{i,j+1/2}=\frac{\rho_{ij}+\rho_{i,j+1}}{2},\\
	{(u_x)}_{i+1/2, j}&=\frac{{(u_x)}_{ij}+{(u_x)}_{i+1,j}}{2},\qquad  {(u_y)}_{i,j+1/2}=\frac{{(u_y)}_{ij}+{(u_y)}_{i,j+1}}{2},
	\end{align*}
	where the macroscopic velocity is defined by
	\begin{align}\label{eq:velocitydiscrete}
	{(u_x)}_{ij}=\frac{1}{\Delta x\Delta y}\sum_{k,l} \rho_{kl} {(F_x)}_{ij}^{kl},\qquad {(u_y)}_{ij}=\frac{1}{\Delta x\Delta y}\sum_{k,l} \rho_{kl} {(F_y)}_{ij}^{kl}
	\end{align}
	with 
	\begin{align*}
	{(F_x)}_{ij}^{kl}&=\iint\limits_{C_{kl}}\bl \iint\limits_{C_{ij}} F_x(x-x',y-y',T(x,y)) \di (x, y) \br\di (x', y'),\\
	{(F_y)}_{ij}^{kl}&=\iint\limits_{C_{kl}}\bl \iint\limits_{C_{ij}} F_y(x-x',y-y',T(x,y)) \di (x, y) \br\di (x', y')
	\end{align*}
	for the components $F_x,F_y$ of $F$ with $F=(F_x,F_y)$. 
	A change of variable also yields
	\begin{align*}
	{(u_x)}_{i+1/2,j}=\frac{1}{\Delta x\Delta y}\sum_{k,l} \rho_{k+1/2,l} {(F_x)}_{ij}^{kl},\qquad {(u_y)}_{i,j+1/2}=\frac{1}{\Delta x\Delta y}\sum_{k,l} \rho_{k,l+1/2} {(F_y)}_{ij}^{kl}.
	\end{align*}
	Note that ${(F_x)}_{ij}^{kl}$ and ${(F_y)}_{ij}^{kl}$ can be determined explicitly in the numerical simulations instead of evaluating the integrals, and can also be precomputed for making the computation of the discretised velocity fields more efficient. 
	Further note that the last line of the numerical scheme \eqref{eq:numericalscheme} can be regarded as a discretisation of the nonlinear diffusion $\delta \nabla \cdot (\rho \nabla \rho)=\frac{\delta}{2}(\partial_x^2 \rho^2 +\partial_y^2 \rho^2)$.
	
	\subsection{Properties of the scheme: conservation of mass, positivity, convergence}\label{sec:convergence} 
	In \cite{Carrillo2016304}, the convergence of a finite volume method is shown for general measure solutions of the (isotropic) aggregation equation with mildly singular potentials. In this section, we establish a CFL condition for the numerical scheme \eqref{eq:numericalscheme} for the anisotropic aggregation equation \eqref{eq:macroscopiceqnonlin} and prove its weak convergence. 
	
	\begin{lemma}\label{lem:rhobound}
		Let $\rho^{in}\in\mathcal{P}_2(\R^2)$ and define $\rho_{ij}^0$ by \eqref{eq:initialconddiscrete}. 
		The conservation of mass is satisfied for all $n$, i.e.
		\begin{align*}
		\sum_{i,j\in \Z} \rho_{ij}^n \Delta x \Delta y =\sum_{i,j\in\Z} \rho_{ij}^0\Delta x \Delta y =1.
		\end{align*}
		For spatially homogeneous tensor fields, conservation of the centre of mass also holds, i.e.
		\begin{align*}
		\sum_{i,j\in \Z} x_i\rho_{ij}^n =\sum_{i,j\in\Z}x_i \rho_{ij}^0, \qquad\sum_{i,j\in \Z} y_i\rho_{ij}^n =\sum_{i,j\in\Z}y_i \rho_{ij}^0.
		\end{align*} 
	\end{lemma}	
	\begin{proof}
		The  conservation of mass is directly obtained by summing over $i$ and $j$ in \eqref{eq:numericalscheme}, and noting that $\sum_{i,j\in\Z} \rho_{ij}^0\Delta x \Delta y =1$. 
		The conservation of the centre of mass follows from a discrete integration by parts and the fact that ${(F_x)}_{ij}^{kl}=-{(F_x)}_{kl}^{ij}$ for spatially homogeneous tensor fields. 
	\end{proof}	
	
	For proving the convergence of the numerical scheme, a CFL condition is required:
	\begin{lemma}\label{lem:schemeproperties}	
		Let $\rho^{in}\in\mathcal{P}_2(\R^2)$ and define $\rho_{ij}^0$ by \eqref{eq:initialconddiscrete}. Suppose that the force $F$ is  bounded by $f$ and, given the spatial discretisation $\Delta x, \Delta y$, assume that the $n$th time step $\Delta t^n$ satisfies
		\begin{align}\label{eq:CFL}
		\bl 2f\bl \frac{1}{\Delta x}+\frac{1}{\Delta y}\br + \delta r_n \bl \frac{1}{(\Delta x)^2}+\frac{1}{(\Delta y)^2}\br \br\Delta t^n\leq 1
		\end{align}
		where
		\begin{align*}
		r_n=\sup_{ij} \rho_{ij}^n.
		\end{align*}
		Then the sequences defined in \eqref{eq:numericalscheme}--\eqref{eq:velocitydiscrete} satisfy
		\begin{align*}
		\rho_{ij}^n\geq 0,\qquad |{(u_x)}_{ij}^{n}|\leq f,\qquad |{(u_y)}_{ij}^{n}|\leq f,
		\end{align*}
		for all $i$, $j$ and $n$. 		In particular, there exists a constant $r>0$ such that 
		\begin{align}\label{eq:rhobound}
		\sup_{n,i,j} \rho_{ij}^n\leq r.
		\end{align} 
	\end{lemma}
	\begin{proof}
		By the definition of the velocity \eqref{eq:velocitydiscrete} and the uniform bound $f$ of the force $F$ we obtain 
		\begin{align}\label{eq:velocityboundhelp}
		|{(u_x)}_{ij}^{n}|\leq \Delta x\Delta y f\sum_{k,l} \rho_{kl}^n=f, \qquad |{(u_y)}_{ij}^n|\leq f
		\end{align}
		for all $i,j,n$.
		
		For proving the nonnegativity of the scheme \eqref{eq:numericalscheme}, note that we can rewrite \eqref{eq:numericalscheme} as
		\begin{align}\label{eq:numericalschemerewritten}
		\begin{split}
		\rho_{ij}^{n+1}&=\rho_{ij}^n\bl 1- \frac{\Delta t^n}{\Delta x} \bl \frac{{(u_x)}^n_{i+1/2,j}-{(u_x)}^n_{i-1/2,j}+2f}{2}\br\right.\\
		&\quad\left.-\frac{\Delta t^n}{\Delta y} \bl \frac{{(u_y)}^n_{i,j+1/2}-{(u_y)}^n_{i,j-1/2}+2f}{2}\br-\frac{\delta \Delta t^n}{(\Delta x)^2}\rho_{ij}^n-\frac{\delta \Delta t^n}{(\Delta y)^2}\rho_{ij}^n \br\\
		&\quad+\rho_{i+1, j}^n\frac{\Delta t^n}{2\Delta x}\bl  f-{(u_x)}^n_{i+1/2, j}\br+\rho_{i-1, j}^n\frac{\Delta t^n}{2\Delta x}\bl  f+{(u_x)}^n_{i-1/2, j}\br\\
		&\quad+\rho_{i, j+1}^n\frac{\Delta t^n}{2\Delta y}\bl  f-{(u_y)}^n_{i, j+1/2}\br+\rho_{i, j-1}^n\frac{\Delta t^n}{2\Delta y}\bl  f+{(u_y)}^n_{i,j-1/2 }\br\\&\quad+\frac{\delta \Delta t^n}{2(\Delta x)^2}\bl(\rho_{i+1,j}^n)^2+(\rho_{i-1,j}^n)^2\br+\frac{\delta \Delta t^n}{2(\Delta y)^2}\bl(\rho_{i,j+1}^n)^2+(\rho_{i,j-1}^n)^2\br.
		\end{split}
		\end{align}
		We show the nonnegativity of $\rho_{ij}^n$ by induction on $n$. For $n\in\N$ given, we  assume that $\rho_{ij}^n\geq  0$ for all $i,j\in\Z$. Note that due to condition \eqref{eq:CFL}, all  coefficients in \eqref{eq:numericalschemerewritten} of $\rho_{ij}^n$, $\rho_{i+1,j}^n$, $\rho_{i-1,j}^n$, $\rho_{i,j+1}^n$ and $\rho_{i,j-1}^n$ are nonnegative, and the  terms in the last line are also nonnegative. By induction, we deduce $\rho_{ij}^{n+1}\geq 0$ for all $i,j\in \Z$.  
		
		Since $\rho_{ij}^n\geq 0$, the conservation of mass implies the uniform boundedness of $\rho_{ij}^n$, i.e.\ there exists a constant $r>0$ such that \eqref{eq:rhobound} is satisfied. 
	\end{proof}

	Next, we consider the convergence of the scheme in a weak topology. Let $\mathcal{M}_{loc}(\R^d)$ denote the space of local Borel measures on $\R^d$. For $\rho \in \mathcal{M}_{loc}(\R^d)$, we denote the total variation of $\rho$ by $|\rho|(\R^d)$ and we denote the space of measures in $\mathcal{M}_{loc}(\R^d)$ with finite total variation by $\mathcal{M}_b(\R^d)$. The space of measures $\mathcal{M}_b(\R^d)$ is always endowed with the weak topology $\sigma(\mathcal{M}_b,C_0)$.
	
	Let the characteristic function on some set $[t_n, t_{n+1})\times C_{ij}\subset \R_+\times \R^2$ 
	be denoted by 
	$\chi_{[t_n, t_{n+1})\times C_{ij}}$. For $\Delta=\max\{\Delta x, \Delta y\}$, we define the reconstruction of the discretisation by
	\begin{align*}
	\rho_\Delta(t,x,y)=\sum_{n\in\N}\sum_{i\in\Z}\sum_{j\in\Z} \rho_{ij}^n \chi_{[t_n, t_{n+1})\times C_{ij}}(t,x,y), 
	\end{align*}
	where the boundedness of $\rho_\Delta$ independent of $\Delta$ follows from Lemma~\ref{lem:schemeproperties}.
	Using the definition $u^n_{ij}=({(u_x)}^n_{ij},{(u_y)}^n_{ij})$ in \eqref{eq:velocitydiscrete}, we obtain
	\begin{align*}
	u_{ij}^n=\frac{1}{\Delta x \Delta y}\iint\limits_{C_{ij}} F(\cdot,T(x,y))\ast \rho_\Delta(t_n,\cdot)(x,y)\di (x, y)
	\end{align*}
	and 
	\begin{align*}
	u_\Delta(t,x,y)=\sum_{n\in\N}\sum_{i\in\Z}\sum_{j\in\Z} u_{ij}^n \chi_{[t_n, t_{n+1})\times C_{ij}}(t,x,y).
	\end{align*}
	\begin{theorem}
		Suppose that  the continuous force $F$ is bounded by $f$ and that the tensor field $T$ is continuous. We consider $\rho^{in}\in\mathcal{P}_2(\R^2)$ and define $\rho^0_{ij}$ by \eqref{eq:initialconddiscrete}. 
		Let $S > 0$ be fixed, and suppose that the discretisation in time and space satisfies \eqref{eq:CFL} where $\rho_{ij}^n$ are obtained from \eqref{eq:numericalscheme} for $\rho^0_{ij}$ given. Then, the discretisation $\rho_\Delta$ converges weakly in $\mathcal{M}_b([0, S] \times \R^2)$ towards the solution $\rho$ of \eqref{eq:macroscopiceq} as $\Delta = \max\{\Delta x, \Delta y\}$ and $\delta$  go to $0$ where for each $\Delta$, the sequence of time steps $\{\Delta t^n\}$  satisfies \eqref{eq:CFL}.
	\end{theorem}
	
	\begin{proof}
		Lemma \ref{lem:rhobound} implies the nonnegativity of $\rho_{ij}^n$ provided  condition \eqref{eq:CFL} holds.  By the conservation of  mass, we have that the sequence $\{\rho_\Delta\}_{\Delta>0}$  of nonnegative bounded measures satisfies  $|\rho_\Delta (t)|(\R^2)=1$ for all $t\in[0,S]$. Hence, there exists a subsequence, still denoted by $\{\rho_\Delta\}_{\Delta>0}$, which converges to $\rho$ in the weak topology as $\Delta =\{\Delta x, \Delta y \}$  goes to 0 where  for each $\Delta$, the sequence of time steps $\{\Delta t^n\}$ satisfies \eqref{eq:CFL}. Hence,
		\begin{align*}
		\int_0^S \iint \limits_{\R^2} \phi(t,x,y) \rho_\Delta (t,x,y) \di (x, y) \di t \to 	\int_0^S \iint \limits_{\R^2} \phi(t,x,y) \rho (t,x,y) \di (x, y) \di t
		\end{align*}
		for all $\phi\in C_0([0,S]\times \R^2)$.
		
		Let $\Delta x, \Delta y$ and $\Delta=\max\{\Delta x,\Delta y \}$ be given.
		For $S>0$  given, 	note that  $N_S\in\N_{>0}$ and $\Delta t^{N_s-1}>0$ can be chosen such that $S=t_{N_{S}}=\sum_{n=0}^{N_S-1} \Delta t^n$ and condition \eqref{eq:CFL} are satisfied. We set 
		$\Delta t=\min_n \Delta t^n$ and choose $N$ such that $S=N\Delta t$. 
		Let $\mathcal{D}([0,S]\times \R^2)$ denote the space of smooth, compactly supported test functions on $[0,S]\times \R^2$ and for  $s_n=n\Delta t$  consider
		\begin{align*}
		\phi_{ij}^n=\int_{s_n}^{s_{n+1}}\iint \limits_{C_{ij}} \phi(t,x,y)\di (x, y) \di t. 
		\end{align*} 
		Note that $\rho_\Delta (s_{n+1},x_i,y_j)-\rho_\Delta(s_n,x_i,y_j)\in\{0, \rho_{ij}^{\sigma(n)+1}-\rho_{ij}^{\sigma(n)}\}$ for any $\sigma(n)\in\{0,\ldots, N_S-1\}$. Here, $\sigma(n)$ is an increasing function defined iteratively with $\sigma(0)=0$ and $\sigma(n+1)=\sigma(n)$ if $\rho_\Delta (s_{n+1},x_i,y_j)-\rho_\Delta(s_n,x_i,y_j)=0$ and $\sigma(n+1)=\sigma(n)+1$ if $\rho_\Delta (s_{n+1},x_i,y_j)-\rho_\Delta(s_n,x_i,y_j)=\rho_{ij}^{\sigma(n)+1}-\rho_{ij}^{\sigma(n)}$. We define $\tilde{\phi}_{ij}^{\sigma(n)}=\phi_{ij}^n$ if $\rho_\Delta (s_{n+1},x_i,y_j)-\rho_\Delta(s_n,x_i,y_j)=\rho_{ij}^{\sigma(n)+1}-\rho_{ij}^{\sigma(n)}$. 
		In particular, we have 
		\begin{align*}
		&\frac{1}{\Delta t} \sum_{n=0}^{N-1}\sum_{i,j\in\Z} \bl \rho_\Delta (s_{n+1},x_i,y_j)-\rho_\Delta(s_n,x_i,y_j) \br\phi_{ij}^n=\frac{1}{\Delta t}\sum_{n=0}^{N_S-1}\sum_{i,j\in\Z}  \bl \rho_{ij}^{n+1}-\rho_{ij}^{n} \br\tilde{\phi}_{ij}^n
		\\\quad&=-\sum_{n=0}^{N_S-1}\sum_{i,j\in \Z} \rho_{ij}^n \frac{\tilde{\phi}_{ij}^n-\tilde{\phi}_{ij}^{n-1}}{\Delta t}=-\sum_{n=0}^{N-1}\sum_{i,j\in \Z} \rho_{ij}^{\sigma(n)} \frac{\phi_{ij}^n-\phi_{ij}^{n-1}}{\Delta t}\\\quad&=-\sum_{n=0}^{N-1}\int_{s_n}^{s_{n+1}}\iint \limits_{\R^2}\rho_\Delta (t,x,y)\frac{\phi(t,x,y)-\phi(t-\Delta t,x,y)}{\Delta t}\di (x,  y) \di t\\&\to -\int_{0}^S\iint \limits_{\R^2}\rho (t,x,y)\partial_t\phi(t,x,y)\di (x, y) \di t
		\end{align*}
		as $\Delta $ goes to 0, where the  limit integral follows from $\phi(t, x, y) - \phi(t - \Delta t, x, y) = \partial_t \phi(t, x, y)\Delta t + \mathcal{O}((\Delta t)^2)$, the weak convergence of $\rho_{\Delta}$ to $\rho$ and the boundedness of the measure $\rho_{\Delta}$  with a bound not depending on the mesh. Note that for 
		\begin{align*}
		\phi_{ij}^n=\int_{t_n}^{t_{n+1}}\iint \limits_{C_{ij}} \phi(t,x,y)\di (x, y) \di t 
		\end{align*} 
		we have
		\begin{align*}
		&\sum_{n=0}^{N_S-1}\sum_{i,j\in\Z} \frac{1}{2\Delta x} \bl \rho_\Delta (t_{n},x_{i+1},y_j)-2\rho_\Delta(t_n,x_i,y_j)+\rho_\Delta (t_{n},x_{i-1},y_j) \br\phi_{ij}^n\\&=\int_{0}^S\iint \limits_{\R^2}\rho_\Delta (t,x,y)\frac{\phi(t,x+\Delta x,y)-2\phi(t,x,y)+\phi(t,x-\Delta x,y)}{2\Delta x}\di (x, y) \di t\to 0
		\end{align*}
		as $\Delta $ goes to 0
		since $|\phi(t,x+\Delta x,y)-2\phi(t,x,y)+\phi(t,x-\Delta x,y)|\leq \| \partial_{xx}\phi\|_\infty (\Delta x)^2$.  
		Due to the boundedness of the force $F(\cdot,T(x))$, we can show in a similar way as in  \cite{Carrillo2016304} that
		\begin{align*}
		&\sum_{n=0}^{N_S-1}\sum_{i,j\in\Z} \frac{1}{\Delta x} \bl {(u_x)}^n_{i+1/2,j} \rho^n_{i+1/2j} - {(u_x)}_{i-1/2,j}^n \rho^n_{i-1/2,j}\br\phi_{i,j}^n\\&\to -\int_0^S \iint \limits_{\R^2} \partial_x \phi(t,x,y) (F_x(\cdot,T(x,y))\ast \rho(t,\cdot))(x,y)\rho(t,x,y)\di (x,y) \di t
		\end{align*}
		as  $\Delta $ goes to 0 by the continuity of $F=(F_x,F_y)$ and $T$
		where $F_x$ denotes the first component of the force $F$. Further note that we have
		\begin{align*}
		&\delta \sum_{n=0}^{N_S-1}\sum_{i,j\in\Z} \frac{1}{2(\Delta x)^2}\bl(\rho_{i+1,j}^n)^2-2\bl\rho_{ij}^n\br^2+\bl\rho_{i-1,j}^n\br^2\br	\phi_{ij}^n\\&=	 \delta\sum_{n=0}^{N_S-1}\sum_{i,j\in\Z} \frac{1}{2(\Delta x)^2}\bl\rho_{ij}^n\br^2\bl \phi_{i+1,j}^n-2\phi_{ij}^n+\phi_{i-1,j}^n\br
		\\& \leq\frac{1}{2} \delta  \|\partial_{xx}\phi\|_\infty \int_0^S \iint \limits_{\R^2} (\rho_\Delta (t,x,y))^2 \di (x, y) \di t.
		\end{align*}
		The boundedness of $\rho_\Delta$, independent of $\Delta$,  guarantees that the right-hand side goes to 0 as   $\Delta $ goes to 0.
		
		Multiplying \eqref{eq:numericalscheme} by $\phi_{ij}^n$, summing over $n$, $i$, $j$, and taking the limits $\delta$ and $\Delta $ to 0, we obtain 
		\begin{align*}
		\int_{0}^S\iint \limits_{\R^2} [\partial_t\phi(t,x,y) + \nabla \phi(t,x,y) \cdot (F(\cdot,T(x,y))\ast \rho(t,\cdot))(x,y)]\rho(t,x,y)\di (x, y) \di t=0
		\end{align*}
		in the limit,
		i.e.\ $\rho$ is a solution in  the sense of distributions of the anisotropic aggregation equation \eqref{eq:macroscopiceq}.
	\end{proof}

	\section{Numerical results}\label{sec:numresults}
	
	In this section, we show simulation results for solving the anisotropic aggregation equation with nonlinear diffusion \eqref{eq:macroscopiceqnonlin} numerically using the numerical scheme \eqref{eq:numericalscheme}. 
	For the numerical simulations, we consider the force coefficients $f_s$ and $f_l$ in \eqref{eq:totalforcenew} with $f_s=f_R+\chi f_A$  and $f_l=f_R+ f_A$ as suggested in \cite{During2017}, where $f_R$ and $f_A$ are defined in \eqref{eq:repulsionforcemodel} and \eqref{eq:attractionforcemodel}. To be consistent with the work of Kücken and Champod \cite{Merkel}, we assume that the total force \eqref{eq:totalforcenew} defined via the tensor field $T(x,y):=\chi s(x,y)\otimes s(x,y) +l(x,y)\otimes l(x,y)$ in \eqref{eq:tensorfield} exhibits short-range repulsion and long-range attraction along $l$ and repulsion along $s$. In the following, we consider the force coefficients $f_R$ and $f_A$ with the parameter values in \eqref{eq:parametervaluesRepulsionAttraction}. The computational domain is given by $[-0.5,0.5]^2$ with periodic boundary conditions.
	
	\subsection{Spatially homogeneous tensor fields}
\begin{figure}[htbp]
		\subfloat[$10^{-10}$]{\includegraphics[width=0.24\textwidth]{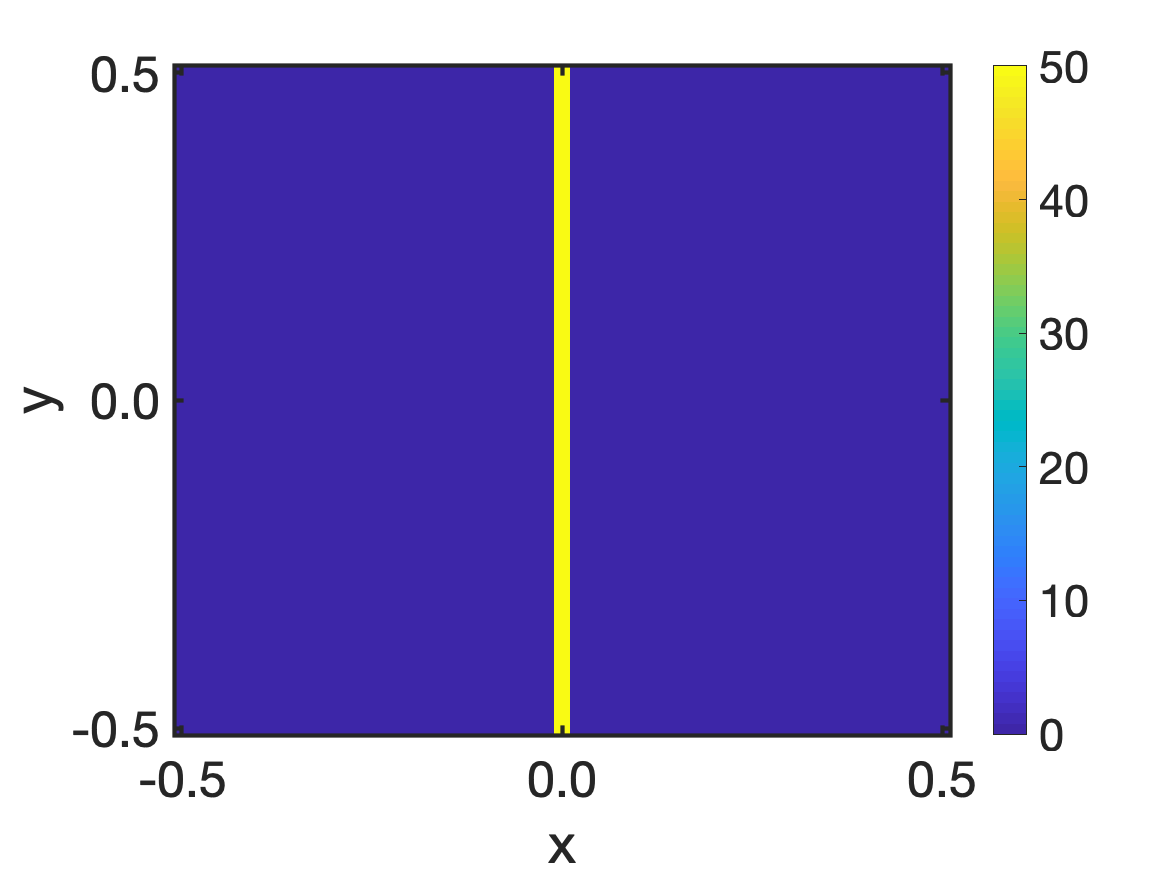}}
		\subfloat[$5\cdot 10^{-8}$]{\includegraphics[width=0.24\textwidth]{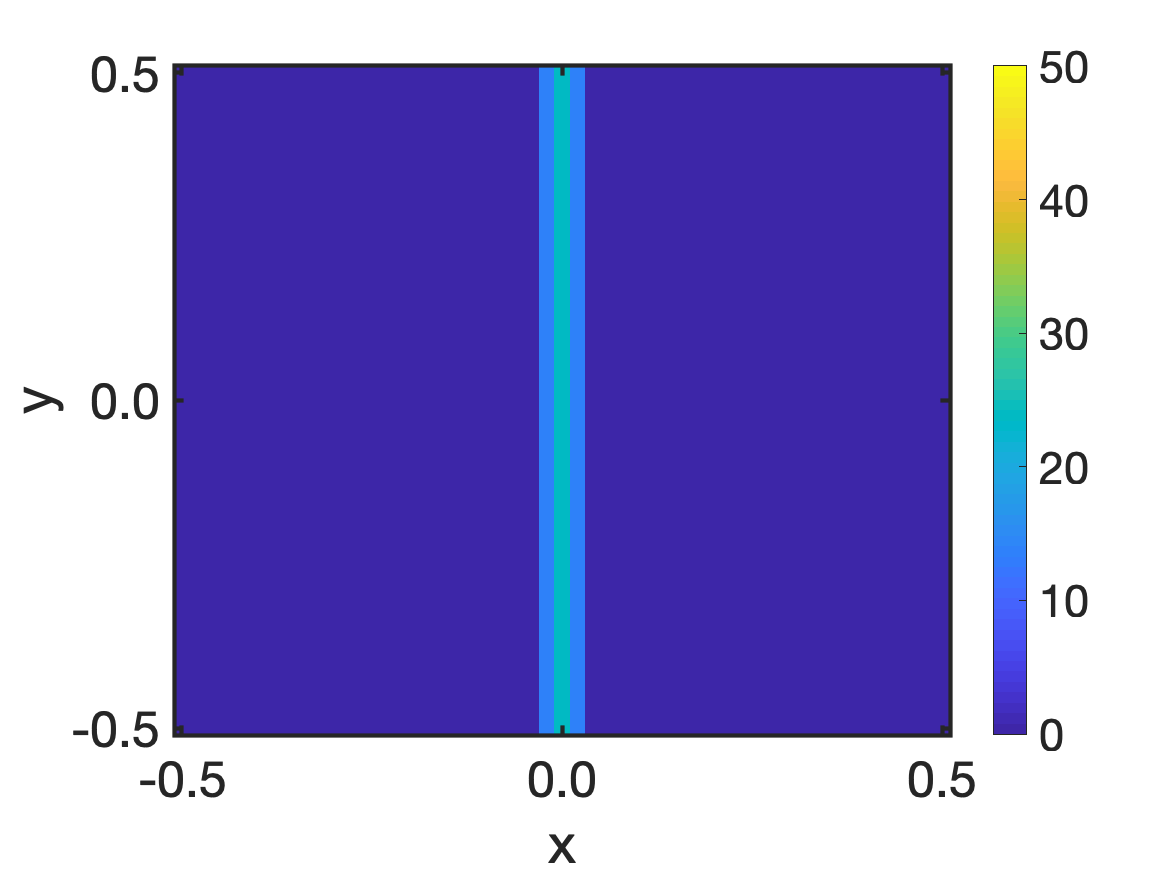}}
		\subfloat[$2\cdot 10^{-7}$]{\includegraphics[width=0.24\textwidth]{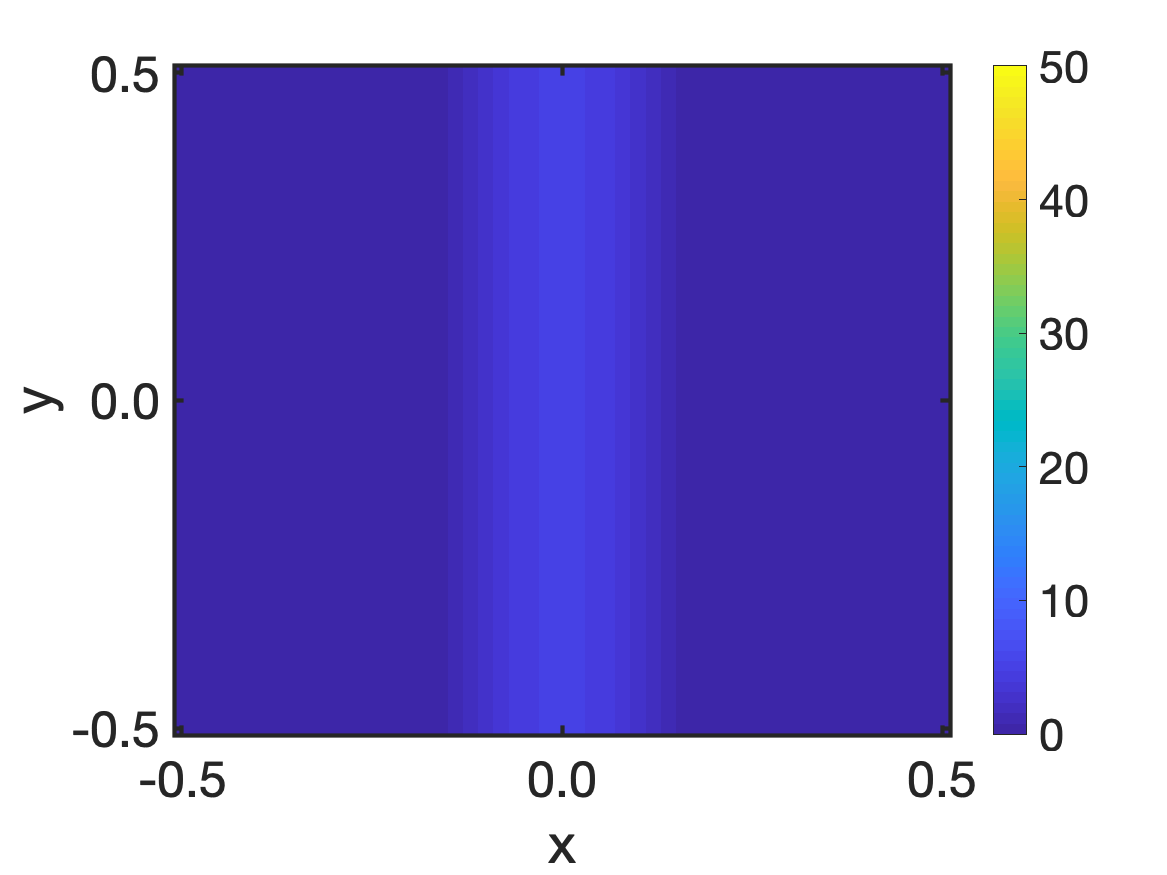}}
		\subfloat[$5\cdot 10^{-7}$]{\includegraphics[width=0.24\textwidth]{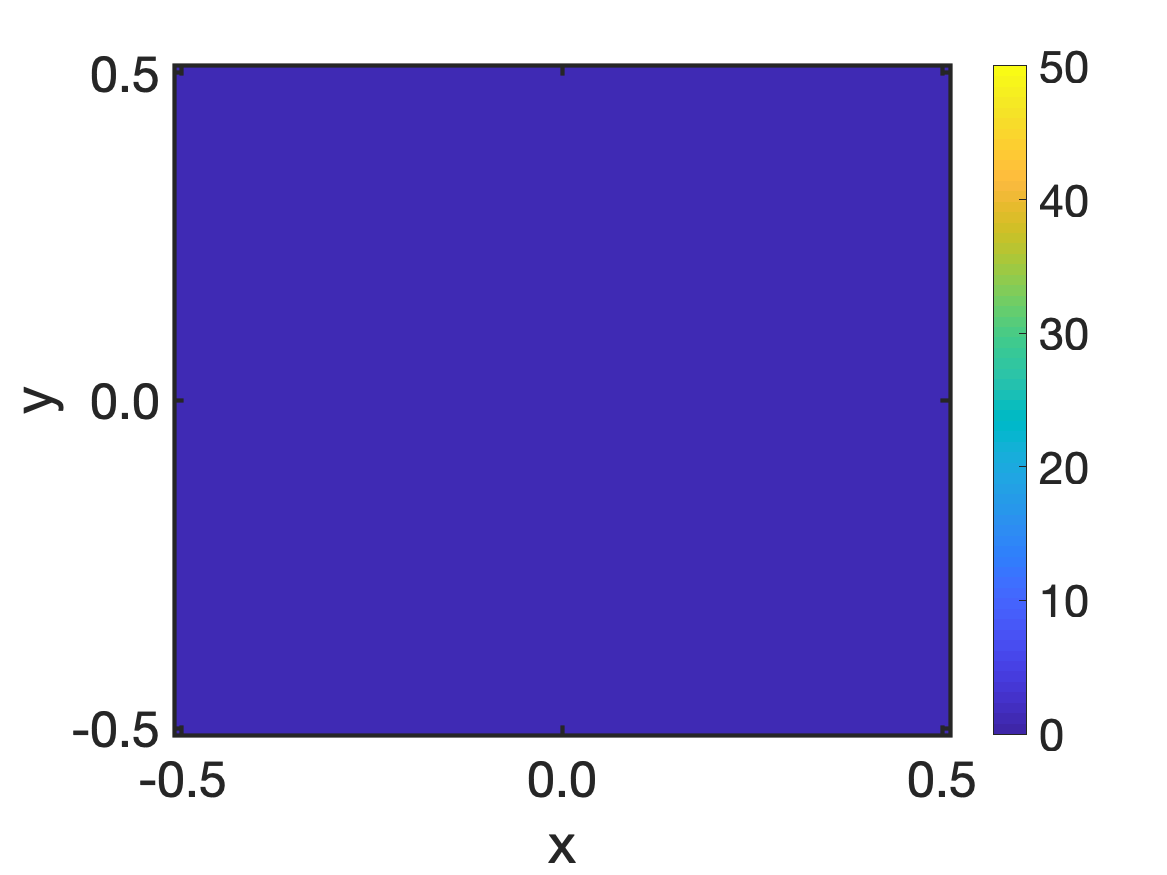}}
		\caption{Stationary solution to the anisotropic interaction equation \eqref{eq:macroscopiceqnonlin}, obtained with the numerical scheme \eqref{eq:numericalscheme} on a grid of size 50 in each spatial direction and different diffusion coefficients $\delta$  for the spatially homogeneous tensor field with $s=(0,1)$ and $l=(1,0)$ and uniformly distributed initial data on a disc on the computational domain $[-0.5,0.5]^2$.}\label{fig:spatiallyhomdiffusion}
	\end{figure}

\begin{figure}[htbp]
		\subfloat[Grid 50]{\includegraphics[width=0.33\textwidth]{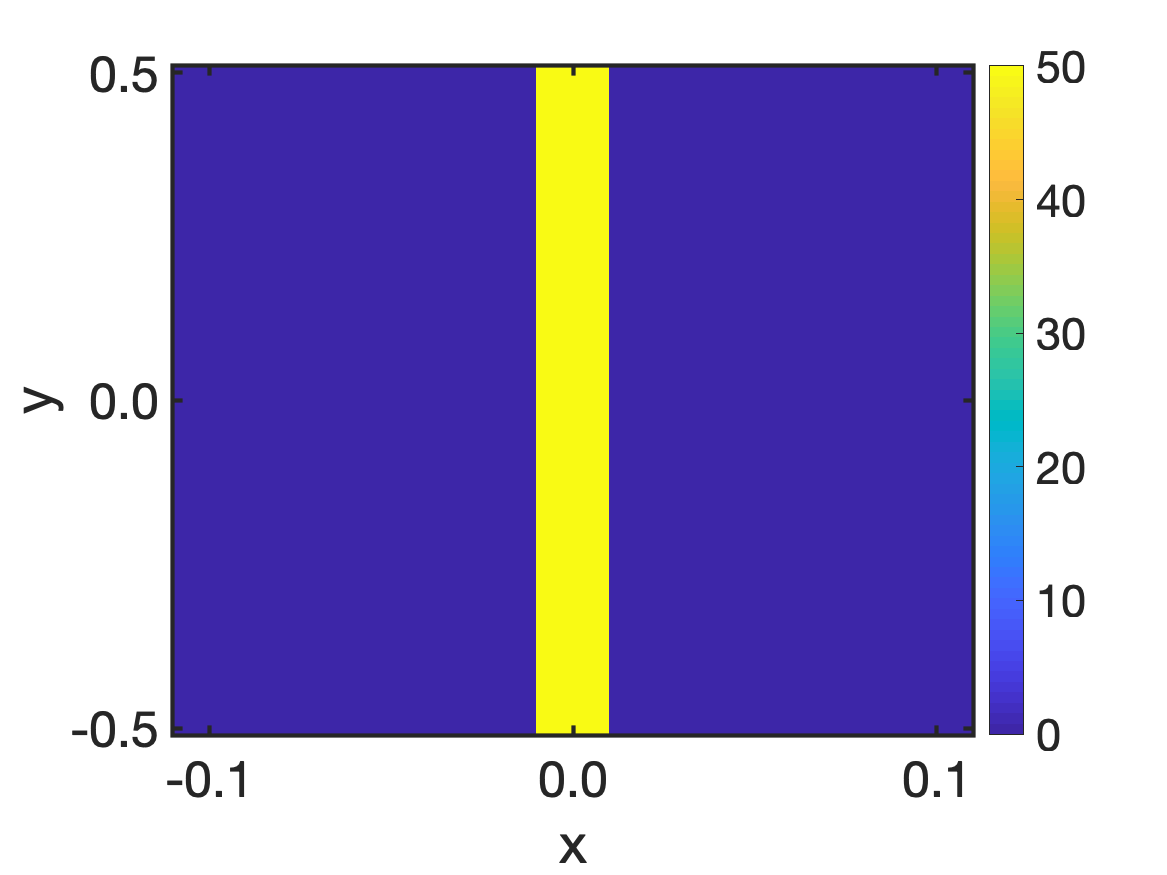}}
		\subfloat[Grid 100]{\includegraphics[width=0.33\textwidth]{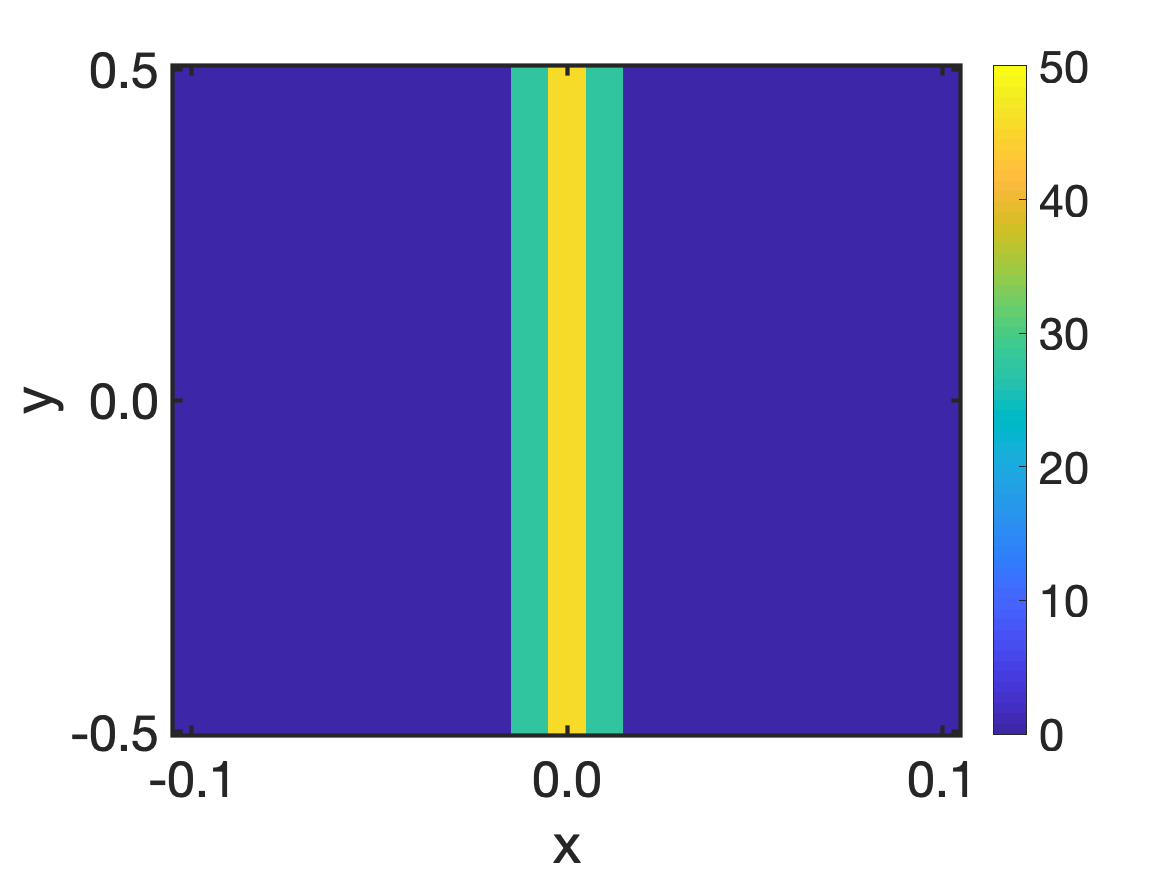}}
		\subfloat[Grid 200]{\includegraphics[width=0.33\textwidth]{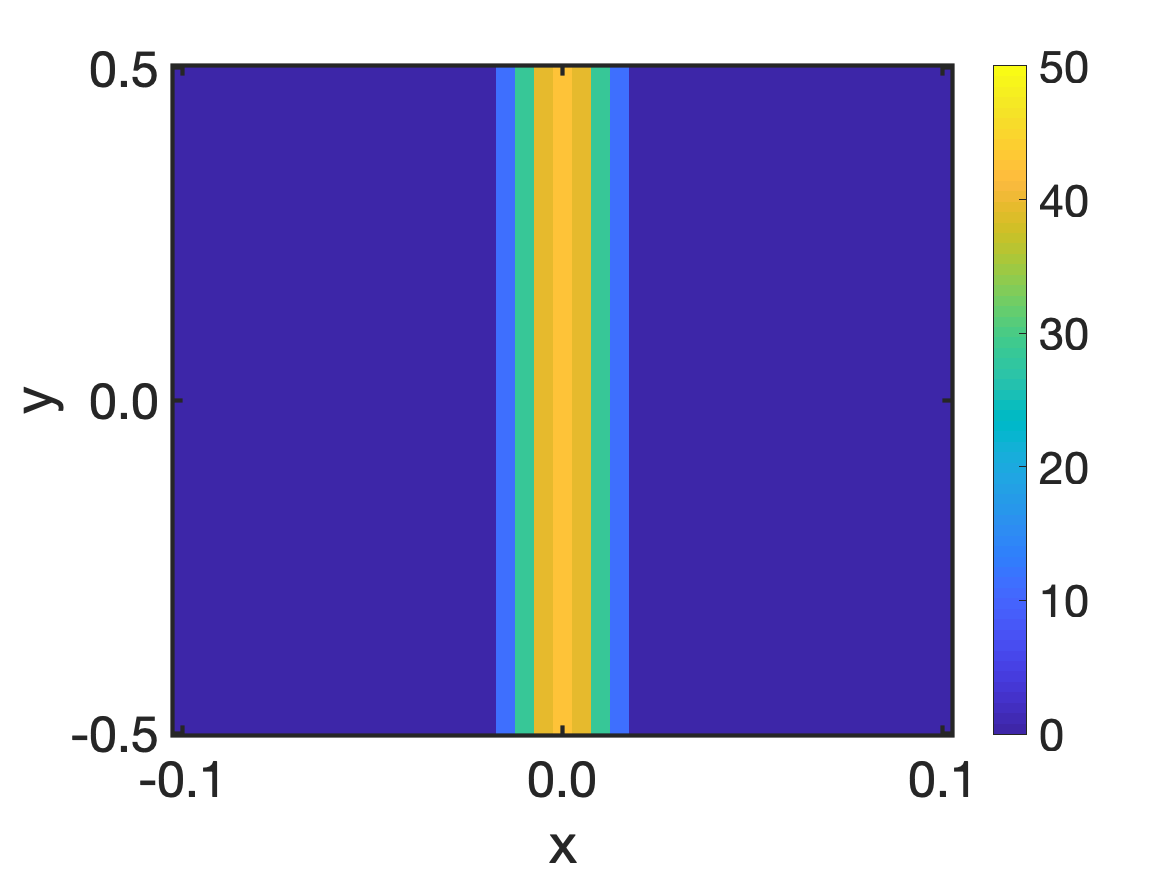}}
		\caption{Stationary solution to the anisotropic interaction equation \eqref{eq:macroscopiceqnonlin}, obtained with the numerical scheme \eqref{eq:numericalscheme} on  grids of sizes 50, 100 and 200 in each spatial direction  for the diffusion coefficient $\delta=10^{-10}$  for the spatially homogeneous tensor field with $s=(0,1)$ and $l=(1,0)$ and uniformly distributed initial data on a disc on the computational domain $[-0.5,0.5]^2$.}\label{fig:spatiallyhomgrid}
	\end{figure}

\begin{figure}[htbp]
		\subfloat[Grid 100]{\includegraphics[width=0.49\textwidth]{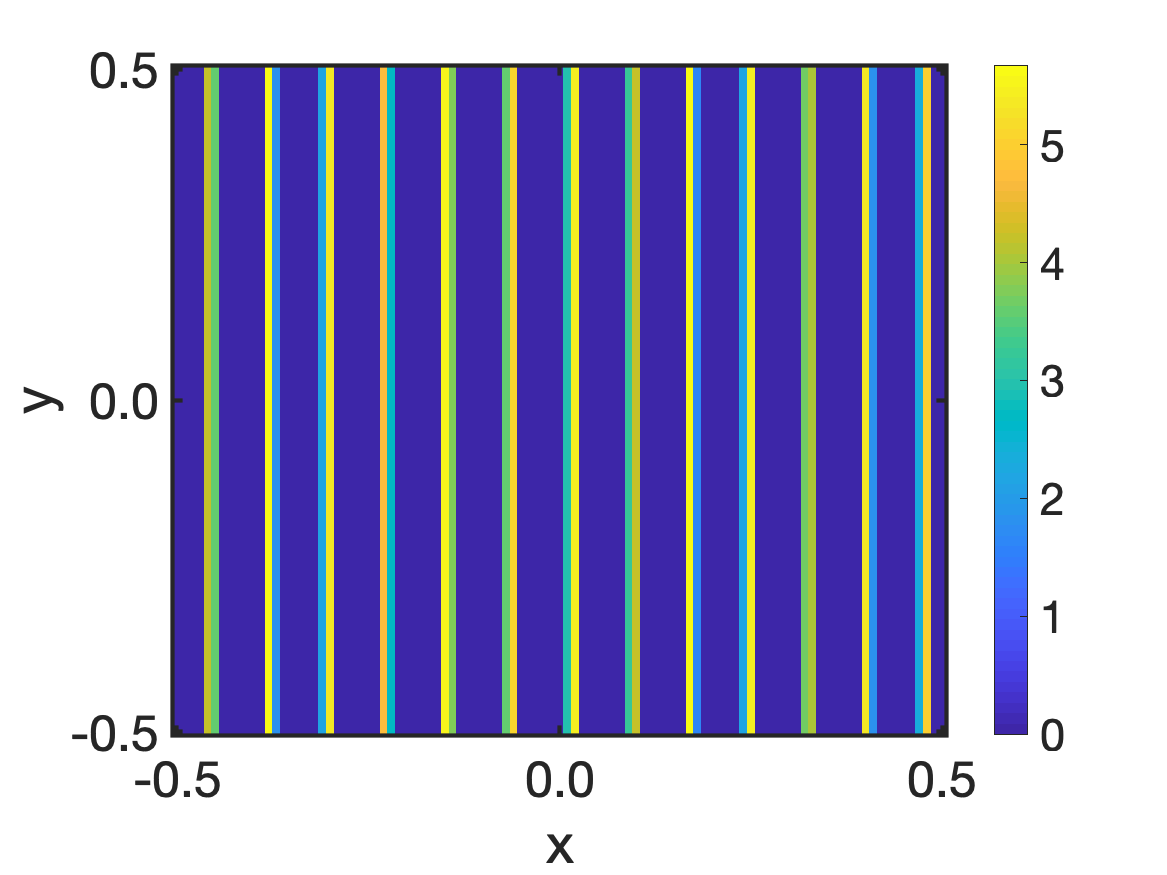}}
		\subfloat[Grid 200]{\includegraphics[width=0.49\textwidth]{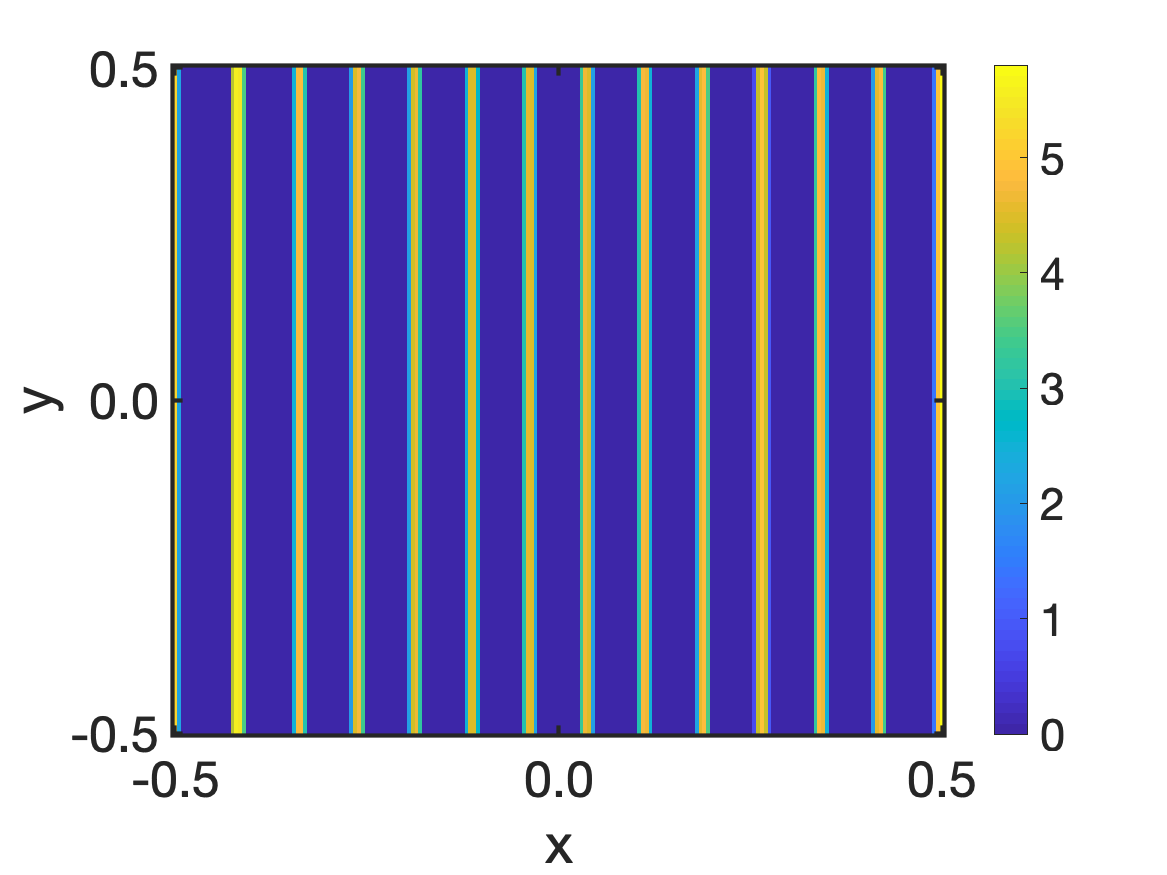}}
		\caption{Stationary solution to the anisotropic interaction equation \eqref{eq:macroscopiceqnonlin}, obtained with the numerical scheme \eqref{eq:numericalscheme} on  grids of sizes  100 and 200 in each spatial direction  for the diffusion coefficient $\delta=10^{-10}$  for the spatially homogeneous tensor field with $s=(0,1)$ and $l=(1,0)$ and uniformly distributed initial data  on the computational domain $[-0.5,0.5]^2$.}\label{fig:spatiallyhomgridvaryinitial}
	\end{figure}

\begin{figure}[htbp]
		\subfloat[$10^{-10}$]{\includegraphics[width=0.24\textwidth]{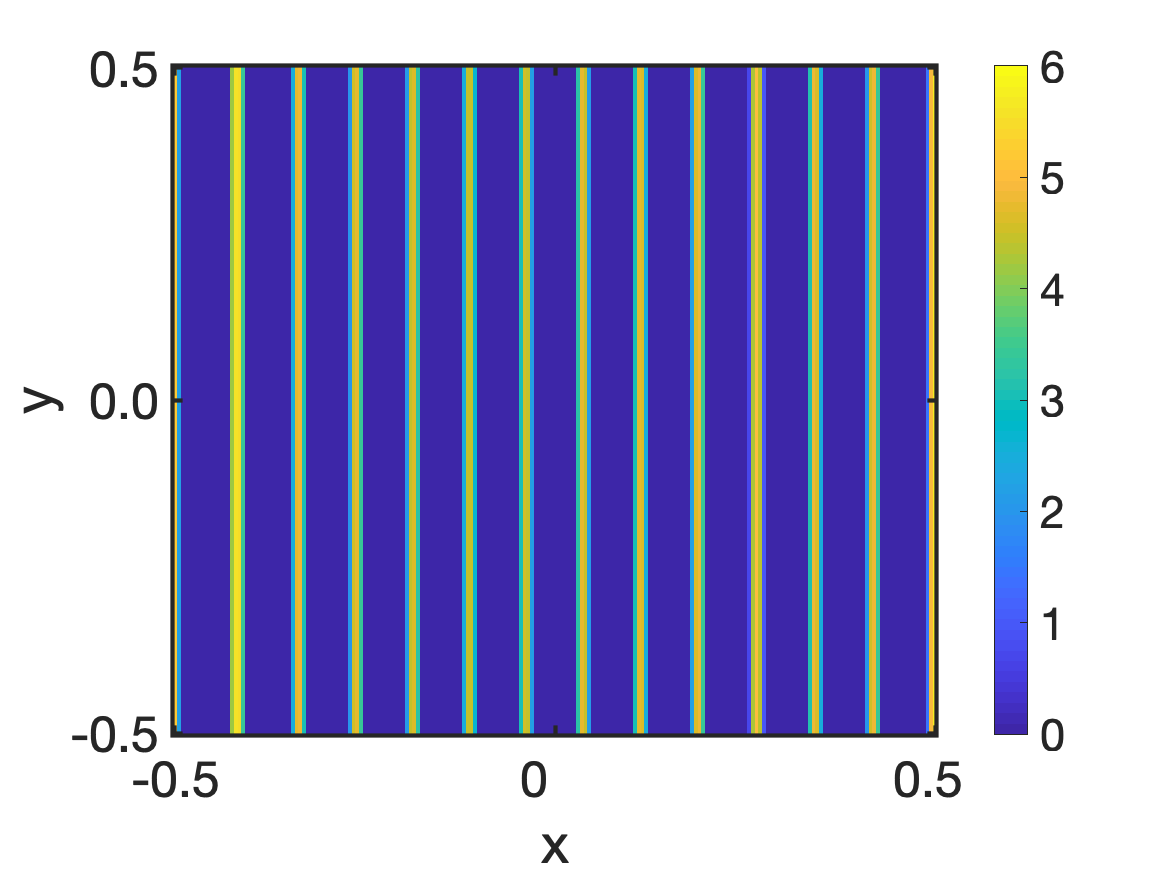}}
		\subfloat[$5\cdot 10^{-10}$]{\includegraphics[width=0.24\textwidth]{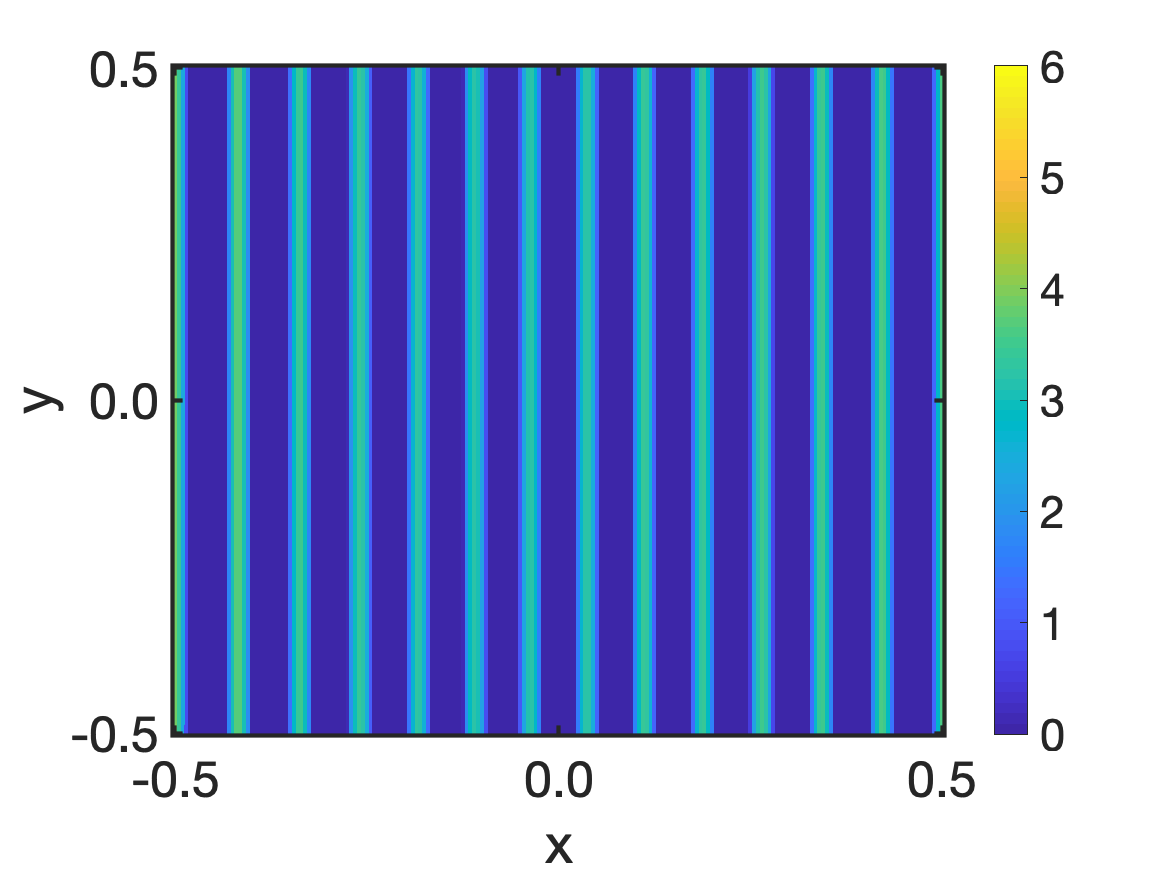}}
		\subfloat[$ 10^{-9}$]{\includegraphics[width=0.24\textwidth]{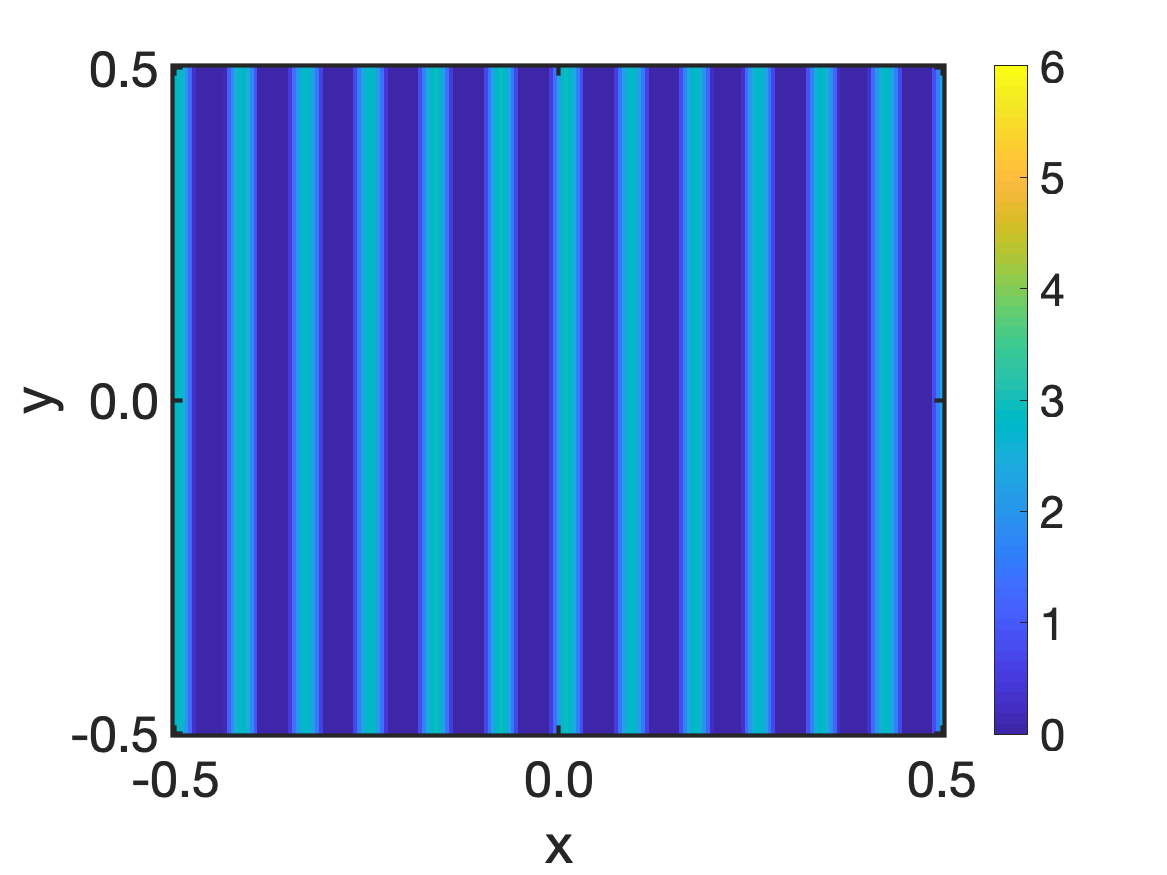}}
		\subfloat[$5\cdot 10^{-9}$]{\includegraphics[width=0.24\textwidth]{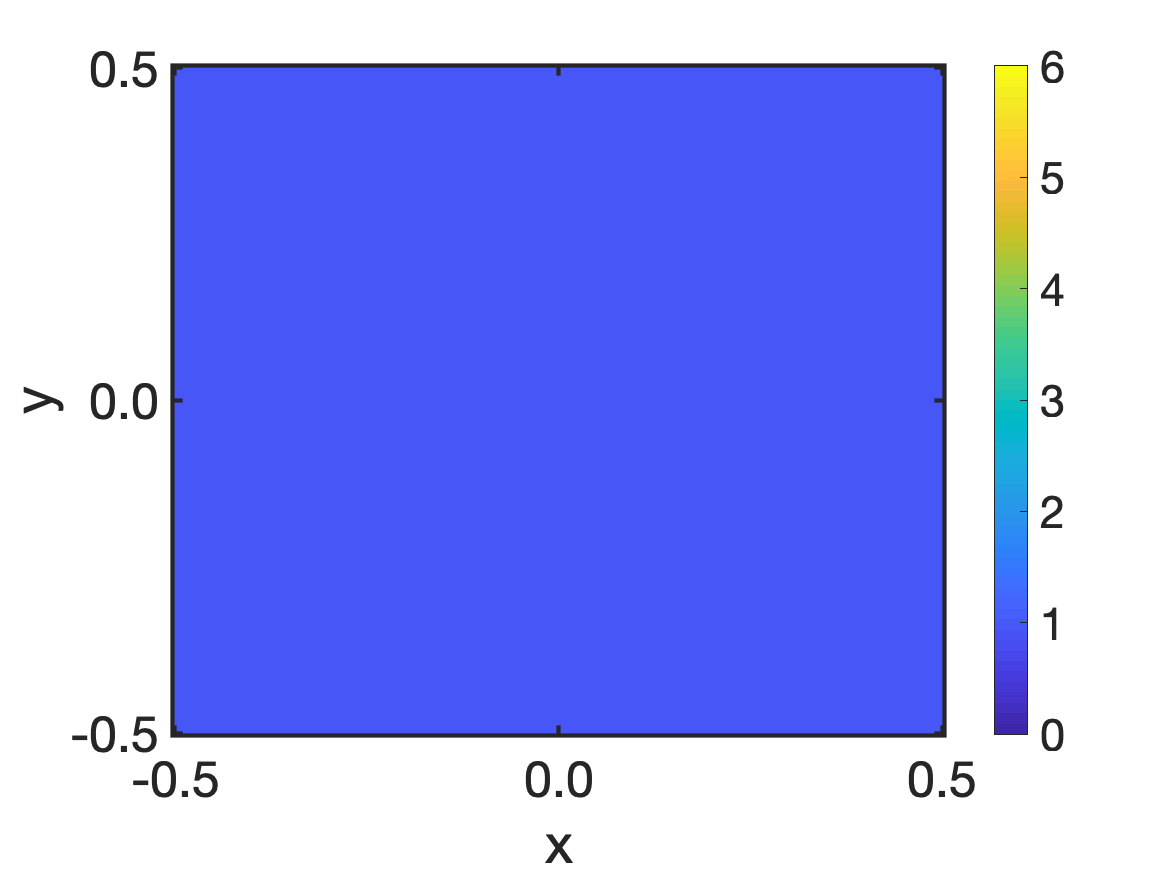}}
		\caption{Stationary solution to the anisotropic interaction equation \eqref{eq:macroscopiceqnonlin}, obtained with the numerical scheme \eqref{eq:numericalscheme} on a grid of size 200 in each spatial direction and different diffusion coefficients $\delta$  for the spatially homogeneous tensor field with $s=(0,1)$ and $l=(1,0)$ and uniformly distributed initial data  on the computational domain $[-0.5,0.5]^2$.}\label{fig:spatiallyhomdiffusionuniformdistributed}
	\end{figure}

	In this section, we show stationary solutions to the anisotropic interaction equation \eqref{eq:macroscopiceqnonlin}, obtained with the numerical scheme \eqref{eq:numericalscheme}, for the spatially homogeneous tensor field $T$ with $s=(0,1)$ and $l=(1,0)$, cf.\ Figures \ref{fig:spatiallyhomdiffusion}--\ref{fig:spatiallyhomdiffusionuniformdistributed}. Note that the stationary solutions for the tensor field $T$ are constant in $y$-direction in all these figures.
	
	The stationary solution to \eqref{eq:macroscopiceqnonlin}, obtained with the numerical scheme \eqref{eq:numericalscheme} for different values of the diffusion coefficient $\delta$,  is shown in Figure \ref{fig:spatiallyhomdiffusion}. Here, we consider uniformly distributed initial data on a disc of radius $R=0.05$ with centre $(0,0)$ on the computational domain $[-0.5,0.5]^2$, where the spatial discretisation is given by a grid of size 50  in each spatial direction, and the time step is chosen according to the CFL condition \eqref{eq:CFL}. Due to the choice of initial data, this leads to a single straight vertical line as stationary solution, provided $\delta$ is chosen sufficiently small. As expected, an increase in $\delta$ leads to the widening of the single straight vertical line which is stable for sufficiently small values of $\delta$. For larger values of $\delta$, e.g.\ $\delta=5\cdot 10^{-7}$, the uniform distribution is obtained as  stationary solution. 
	
	In Figure \ref{fig:spatiallyhomgrid}, we investigate the role of the grid size on the stationary solution by considering  grid sizes of 50, 100 and 200 in each spatial direction for the diffusion parameter $\delta=10^{-10}$ and uniformly distributed initial data on a disc. Clearly, the stationary solution is given by a step function in the $x$-coordinate. Finer grids lead to step functions  with more steps and smaller step heights compared to the grid size of 50 where only one step occurs. 
	
	The stationary solution for grid sizes of 100 and 200 in each spatial direction and uniformly distributed initial data on the computational domain $[-0.5,0.5]^2$ is shown in Figure \ref{fig:spatiallyhomgridvaryinitial}, and is given by equidistant, parallel vertical line patterns. Note that we obtain the same number of parallel lines for the different grid sizes.

\begin{figure}
		\includegraphics[width=0.32\textwidth]{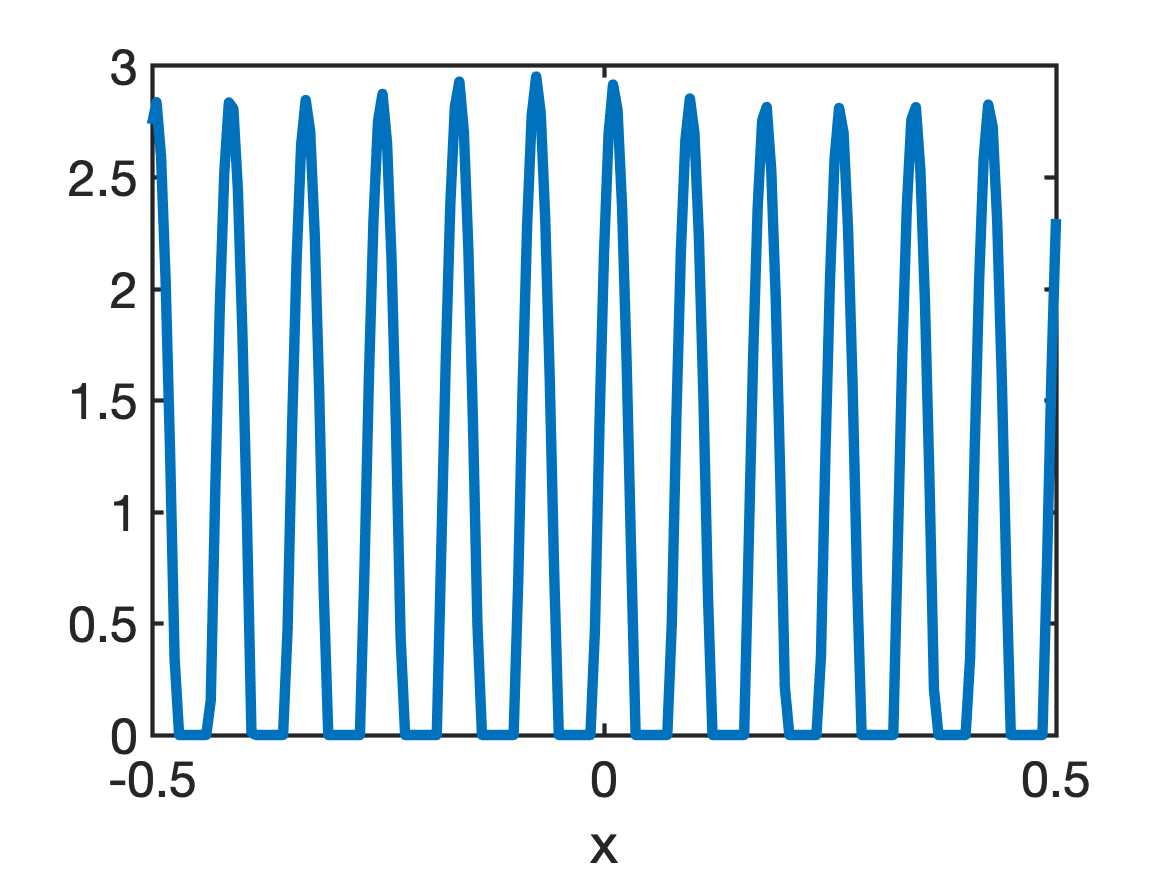}		
		\caption{Cross-section of stationary solution  to the anisotropic interaction equation \eqref{eq:macroscopiceqnonlin}, obtained with the numerical scheme \eqref{eq:numericalscheme} on a grid of size 200 in each spatial direction and  diffusion coefficient $\delta=10^{-9}$  for the spatially homogeneous tensor field with $s=(0,1)$ and $l=(1,0)$ and uniformly distributed initial data  on the computational domain $[-0.5,0.5]^2$.}\label{fig:spatiallyhomdiffusionuniformdistributedcrosssection}
	\end{figure}

	In Figure \ref{fig:spatiallyhomdiffusionuniformdistributed}, we show  the stationary solution for uniformly distributed initial data on the computational domain $[-0.5,0.5]^2$ for different diffusion coefficients $\delta$. Note that as $\delta$ increases, the stable line patterns become wider and this may result in a decrease in the number of parallel lines. If $\delta$ is larger than a certain threshold, e.g.\  $\delta=5\cdot 10^{-9}$, the parallel line patterns are no longer stable and the stationary solution is given by the uniform distribution on the computational domain $[-0.5,0.5]^2$.
	The plot of the cross-section of the stationary solution for diffusion coefficient $\delta=10^{-9}$ is shown in Figure \ref{fig:spatiallyhomdiffusionuniformdistributedcrosssection}. Note that the solution is finite and has no blow-up.

	\subsection{Spatially inhomogeneous tensor fields}
	
	In this section, we consider stationary solutions to the anisotropic interaction equation \eqref{eq:macroscopiceqnonlin}, obtained with the numerical scheme \eqref{eq:numericalscheme},  for different spatially inhomogeneous tensor fields.
	
	This section is motivated by the simulation of fingerprint patterns. The fingerprint development can be modeled in three phases \cite{Merkel}. In the first phase, compressive mechanical stress is created, which has been modeled by K\"{u}cken and Newell \cite{fingerprintformation1,fingerprintformation2}. As proposed in \cite{Merkel}, an easy way to obtain realistic underlying stress fields is to construct a tensor field based on real fingerprint data. For this, we consider real fingerprint images and construct the vector field $s=s(x)$ for all $x\in\Omega$ as the tangents to the given fingerprint lines via extrapolation \cite{During2017,LineSensor}. Given the underlying stress field, the second phase of the fingerprint development consists of the rearrangement of Merkel cells from a random configuration into parallel ridges along the lines of smallest compressive stress. This phase describes the pattern formation, first modeled by K\"{u}cken and Champod  \cite{Merkel}, and its continuum limit is given by \eqref{eq:macroscopiceq}. The construction of $s$ and solving \eqref{eq:macroscopiceq} are motivated by the idea that $s$ denotes the lines of smallest stress and the solution to \eqref{eq:macroscopiceq} aligns along $s$. Finally, primary ridges are induced in the third phase which is not part of this work.
	
\begin{figure}[htbp]
		\subfloat[Original]{\includegraphics[width=0.32\textwidth]{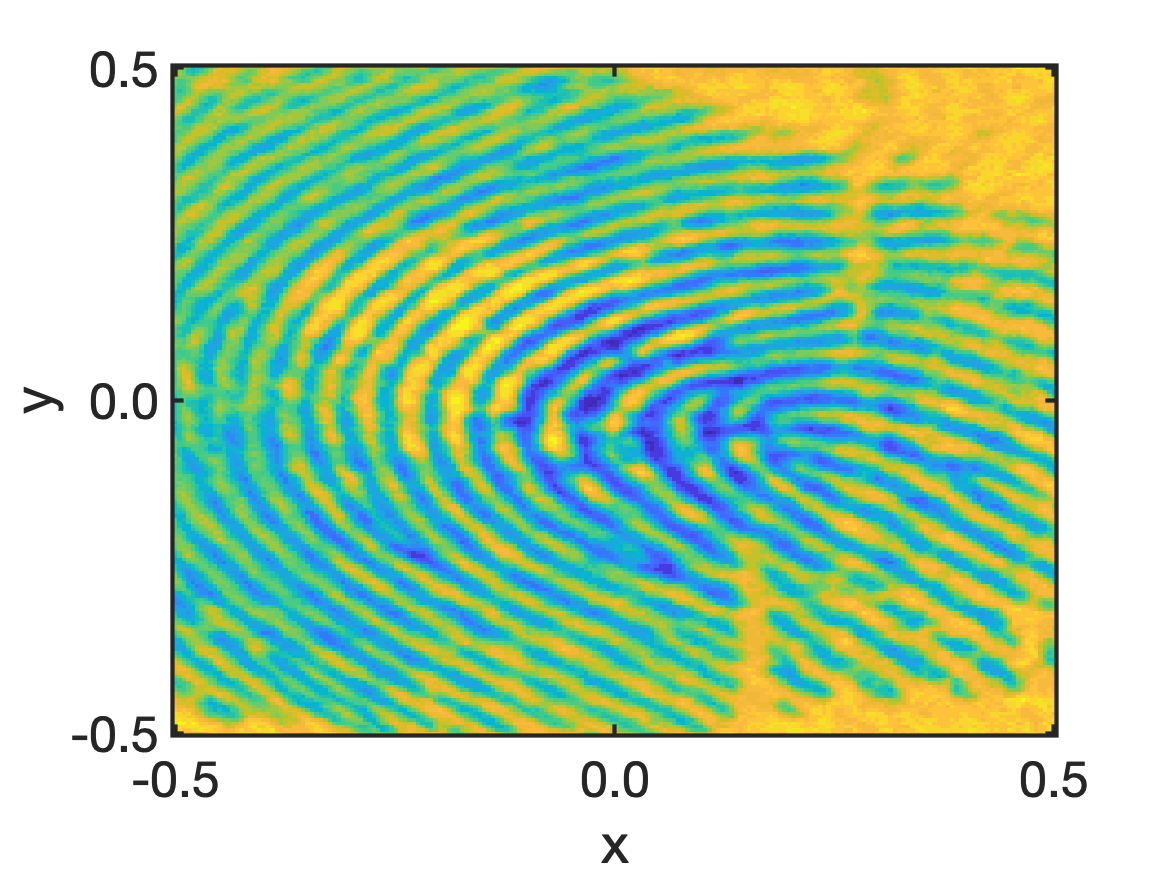}\label{fig:originalpart}}
		\subfloat[s]{\includegraphics[width=0.32\textwidth]{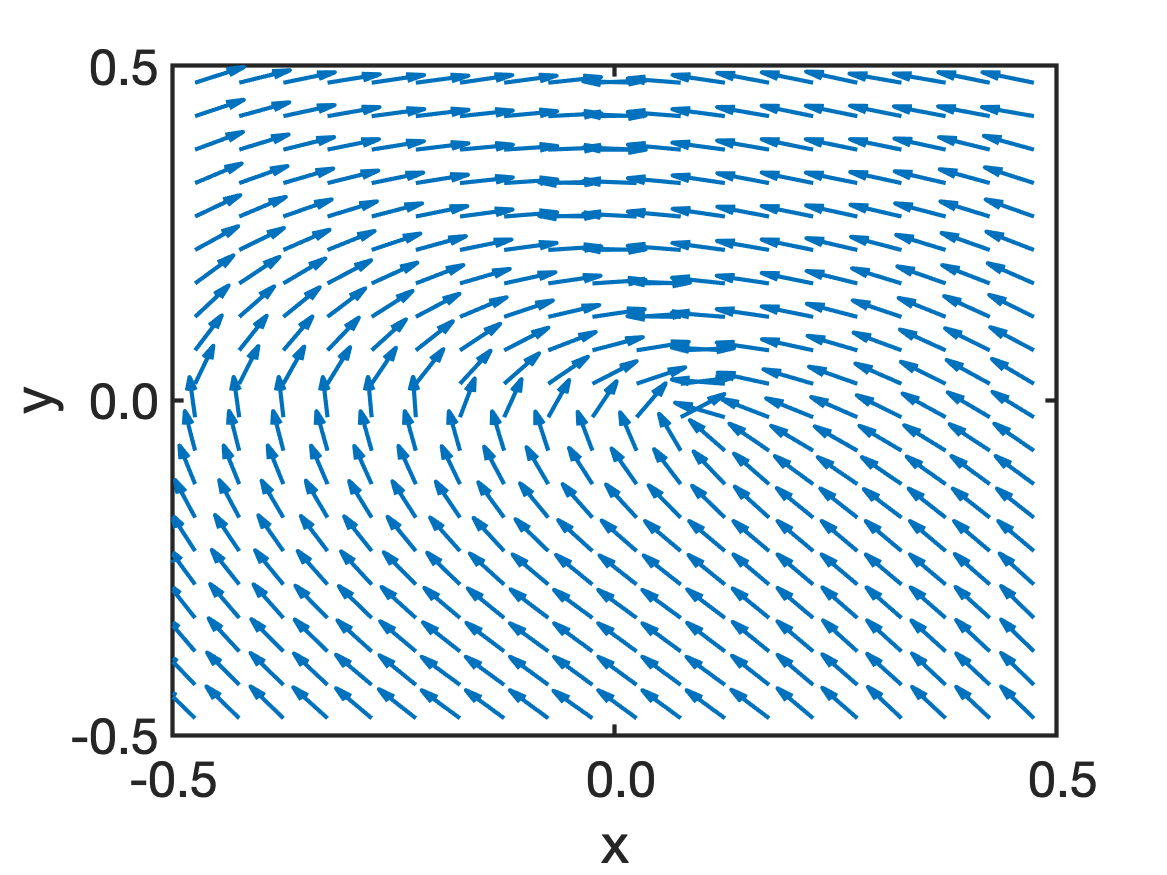}\label{fig:spart}}
		\subfloat[Stationary solution]{\includegraphics[width=0.32\textwidth]{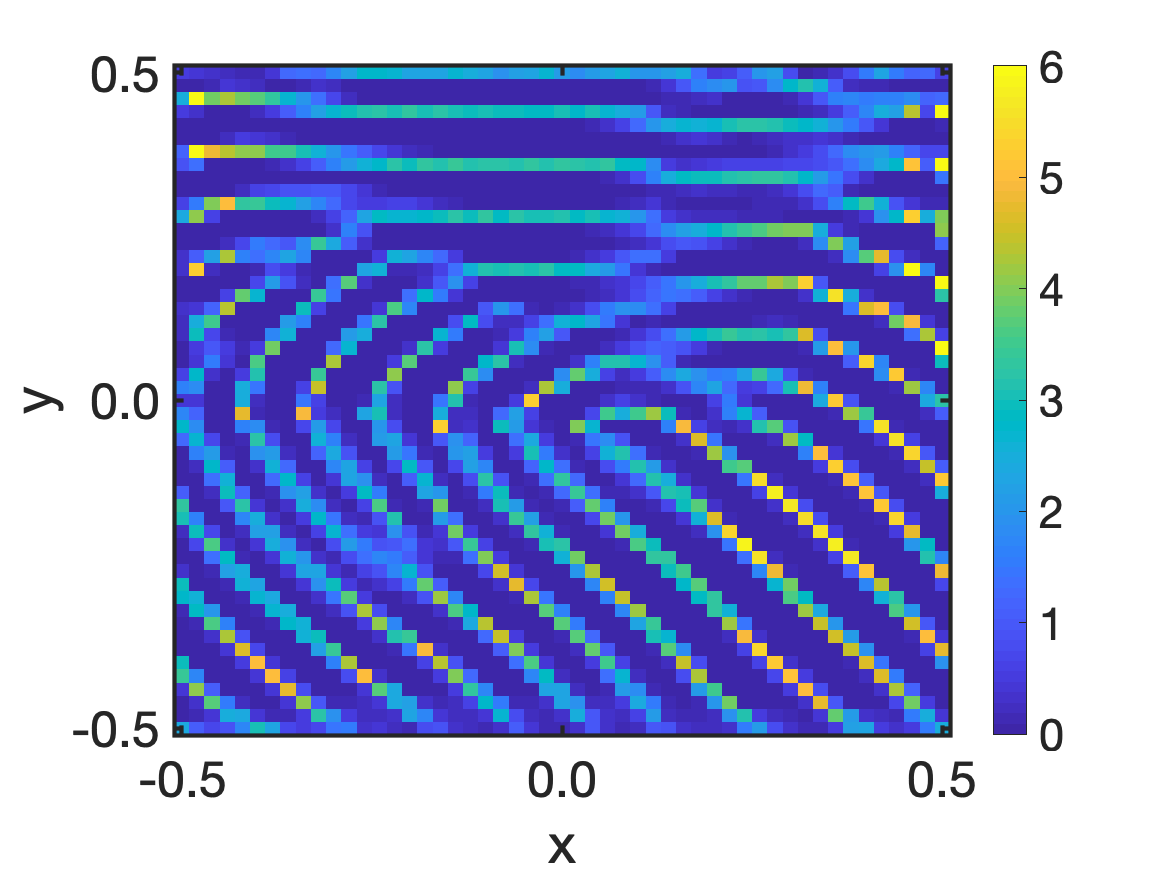}\label{fig:stationarypart}}\\
		\subfloat[Original]{\includegraphics[width=0.32\textwidth]{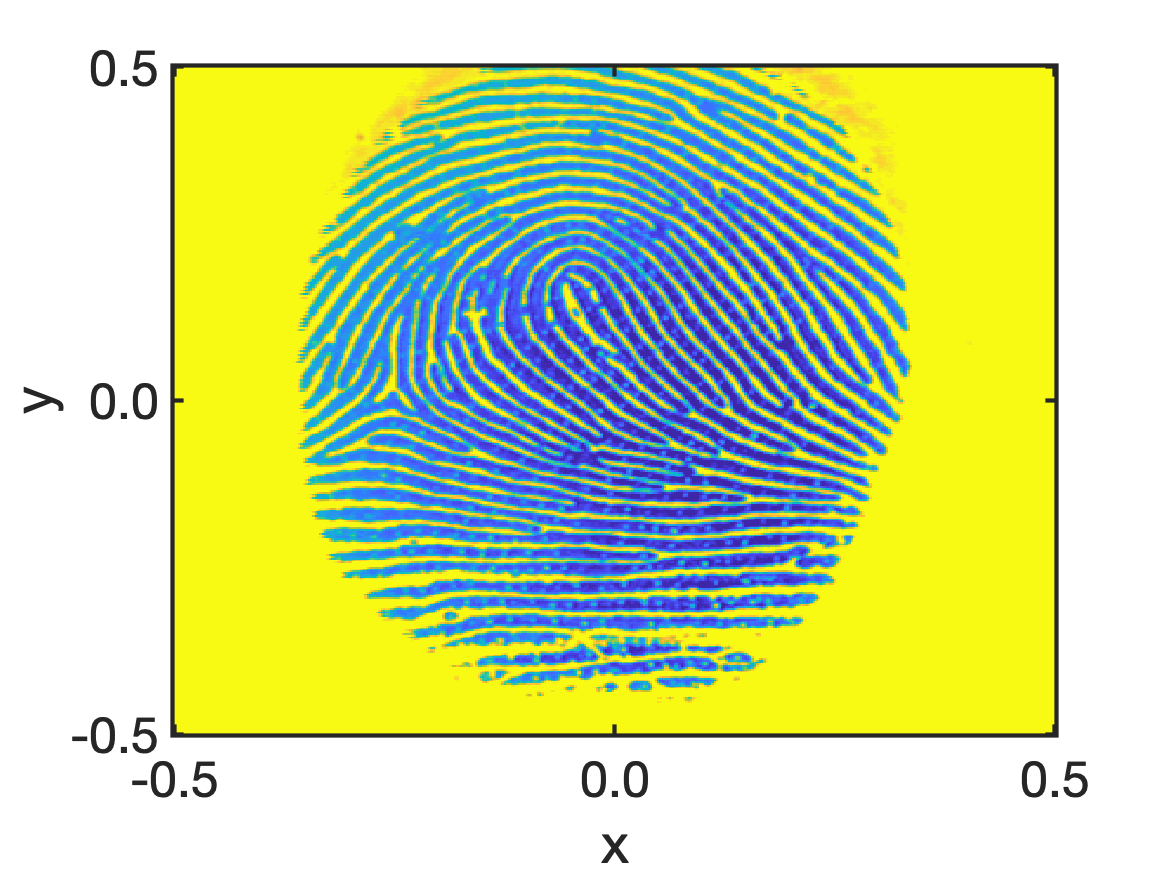}\label{fig:originalfull}}
		\subfloat[s]{\includegraphics[width=0.32\textwidth]{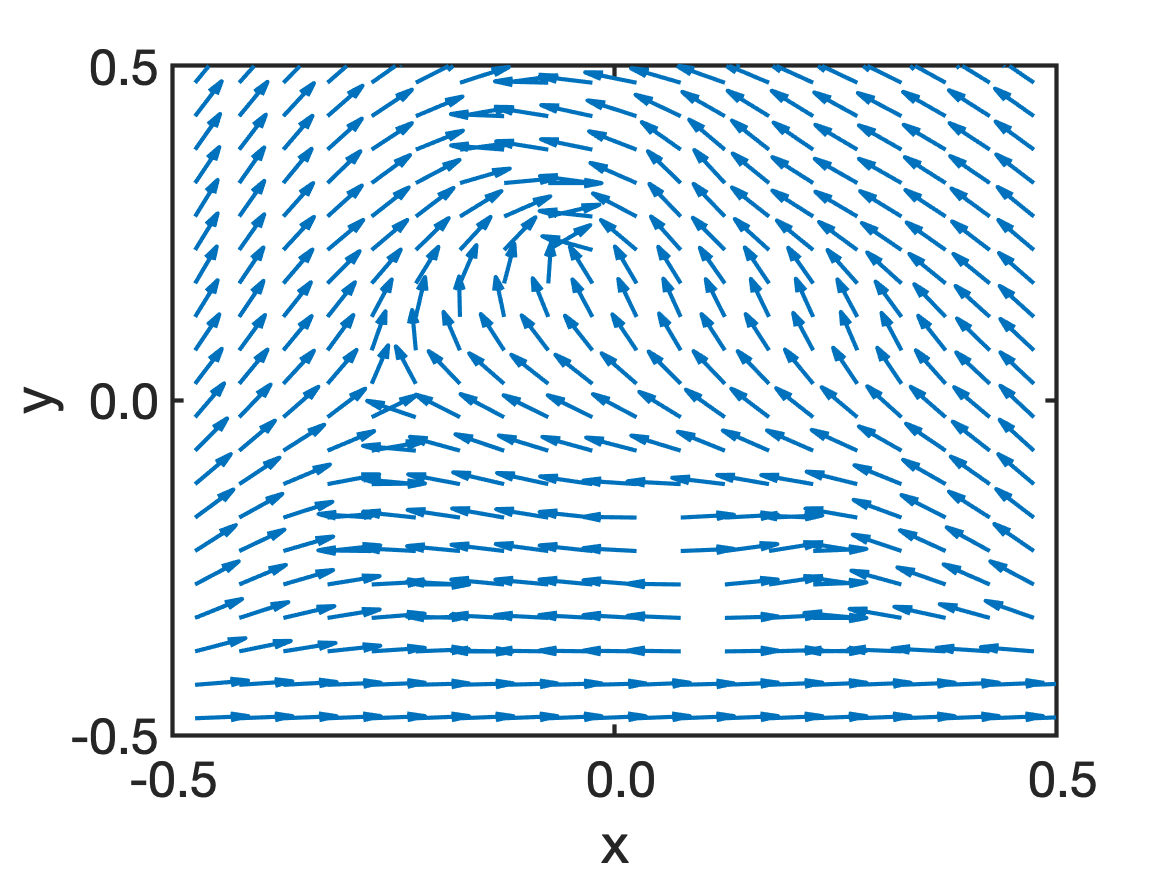}\label{fig:sfull}}
		\subfloat[Stationary solution]{\includegraphics[width=0.32\textwidth]{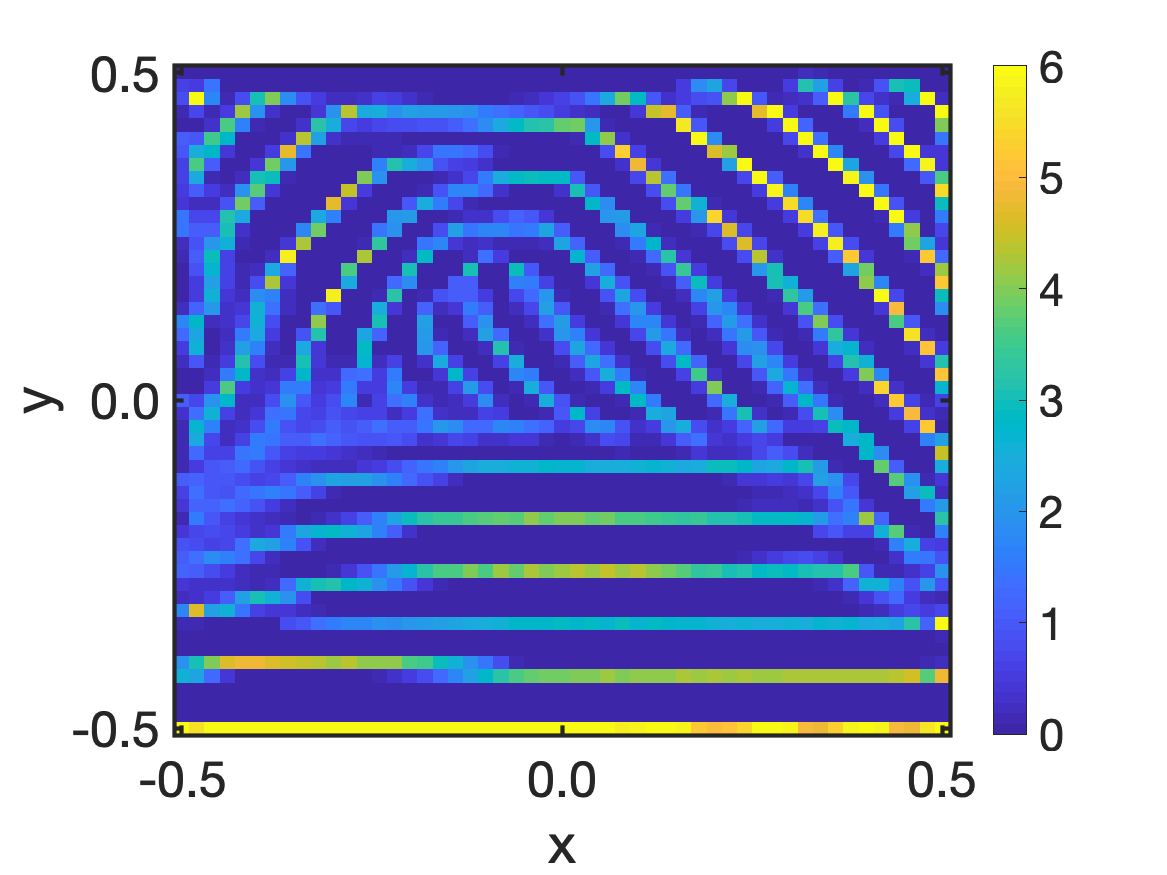}\label{fig:stationaryfull}}
		\caption{Stationary solution to the anisotropic interaction equation \eqref{eq:macroscopiceqnonlin}, obtained with the numerical scheme \eqref{eq:numericalscheme} on a grid of size 50 in each spatial direction and diffusion coefficient $\delta=10^{-10}$  for different spatially inhomogeneous tensor fields from real fingerprint images and uniformly distributed initial data  on the computational domain $[-0.5,0.5]^2$.}\label{fig:stationarytensor}
	\end{figure}

	In Figure \ref{fig:stationarytensor}, we consider fingerprint images in Figures \subref*{fig:originalpart} and \subref*{fig:originalfull}, use these fingerprint images to construct the vector field $s=s(x,y)$ in Figures \subref*{fig:spart} and \subref*{fig:sfull}, and show the resulting stationary solutions for the diffusion coefficient $\delta=10^{-10}$ and  uniformly distributed initial data on a grid of size 50 in each spatial direction in Figures \subref*{fig:stationarypart} and \subref*{fig:stationaryfull}, respectively. For the construction of the tensor field we firstly proceed as in \cite{During2017}, and then we rescale the tensor field appropriately to the given grid size. As desired, the stationary solution in Figures \subref*{fig:stationarypart} and \subref*{fig:stationaryfull} aligns along the vector field $s$ where $s$ is shown in Figures \subref*{fig:spart} and \subref*{fig:sfull}, respectively. Here, the orientation of the stripes is the main feature, while the number of lines depends on the scaling of the forces, also compare Figure \ref{fig:spatiallyInhomdiffusionuniformdistributed}.
	
\begin{figure}[htbp]
		\subfloat[$10^5$]{\includegraphics[width=0.24\textwidth]{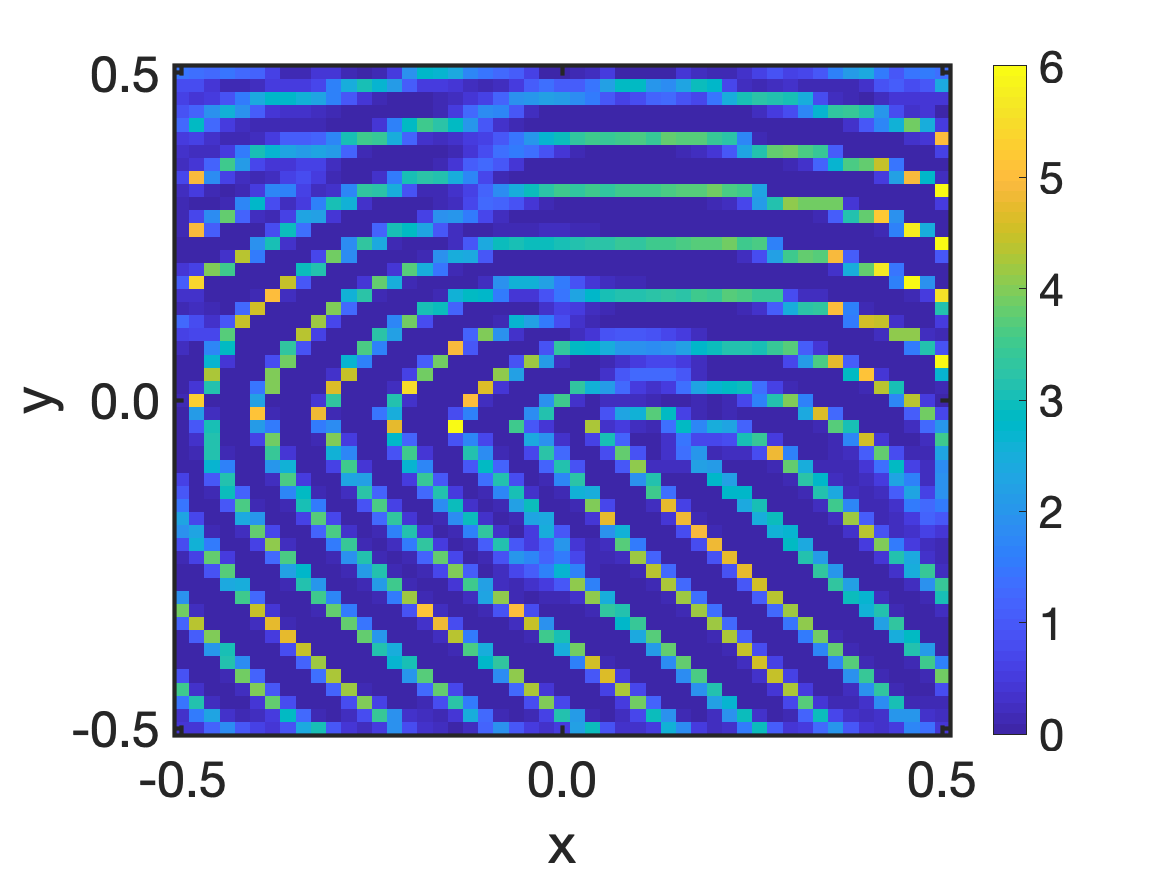}}
		\subfloat[$2\cdot 10^5$]{\includegraphics[width=0.24\textwidth]{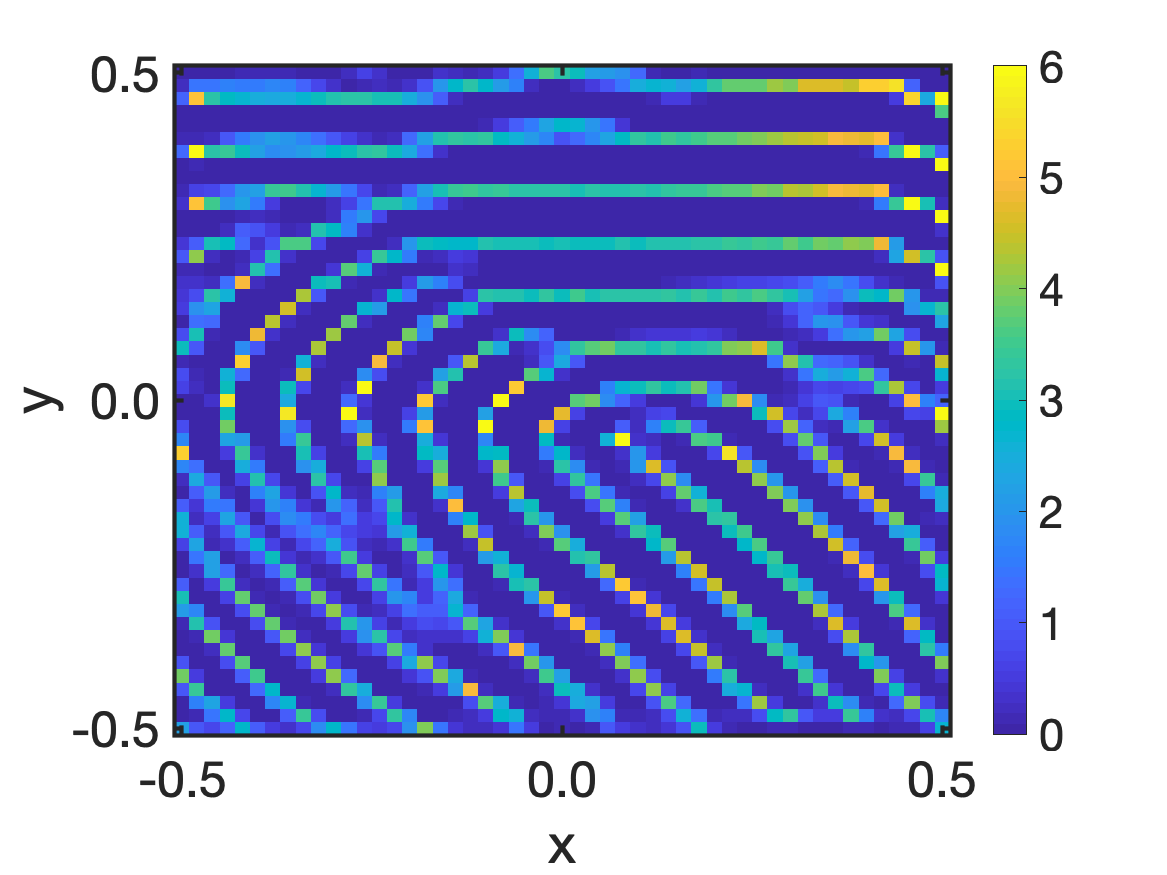}}
		\subfloat[$3\cdot 10^5$]{\includegraphics[width=0.24\textwidth]{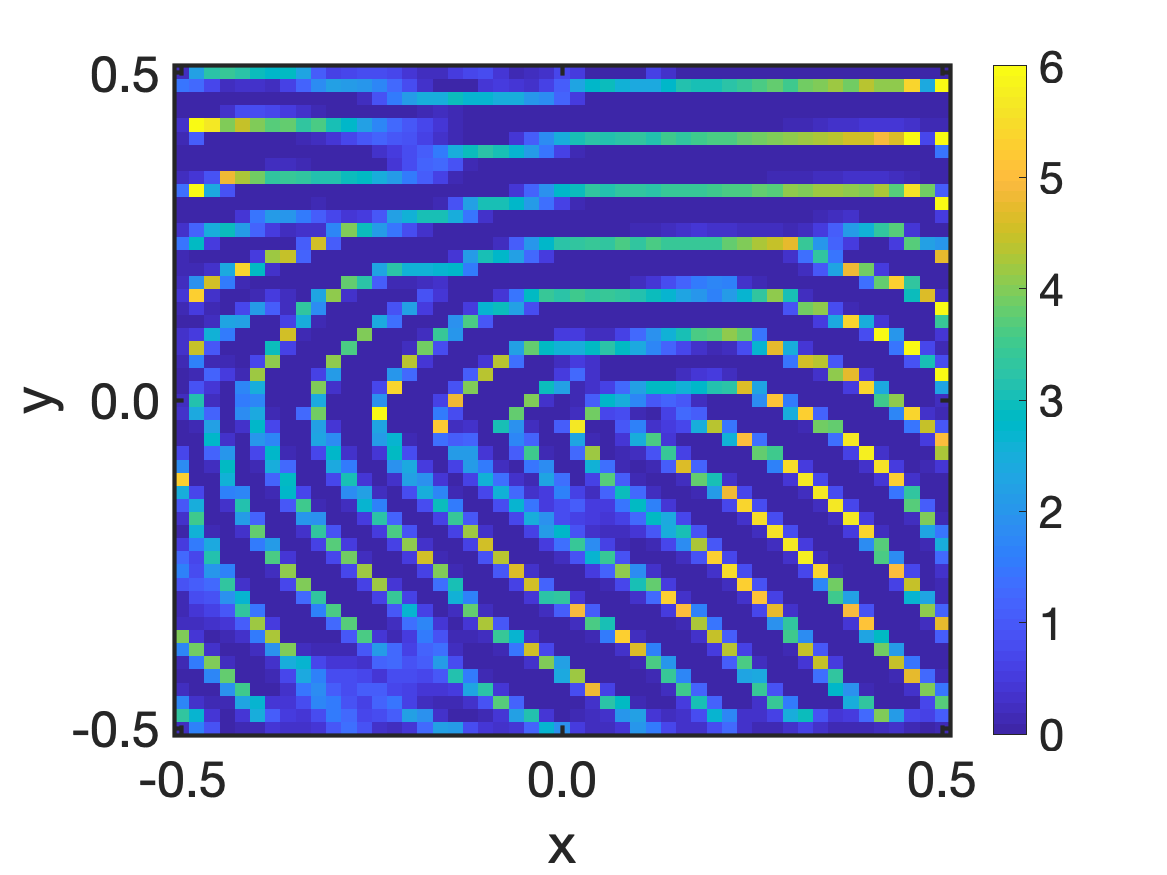}}
		\subfloat[$4\cdot 10^5$]{\includegraphics[width=0.24\textwidth]{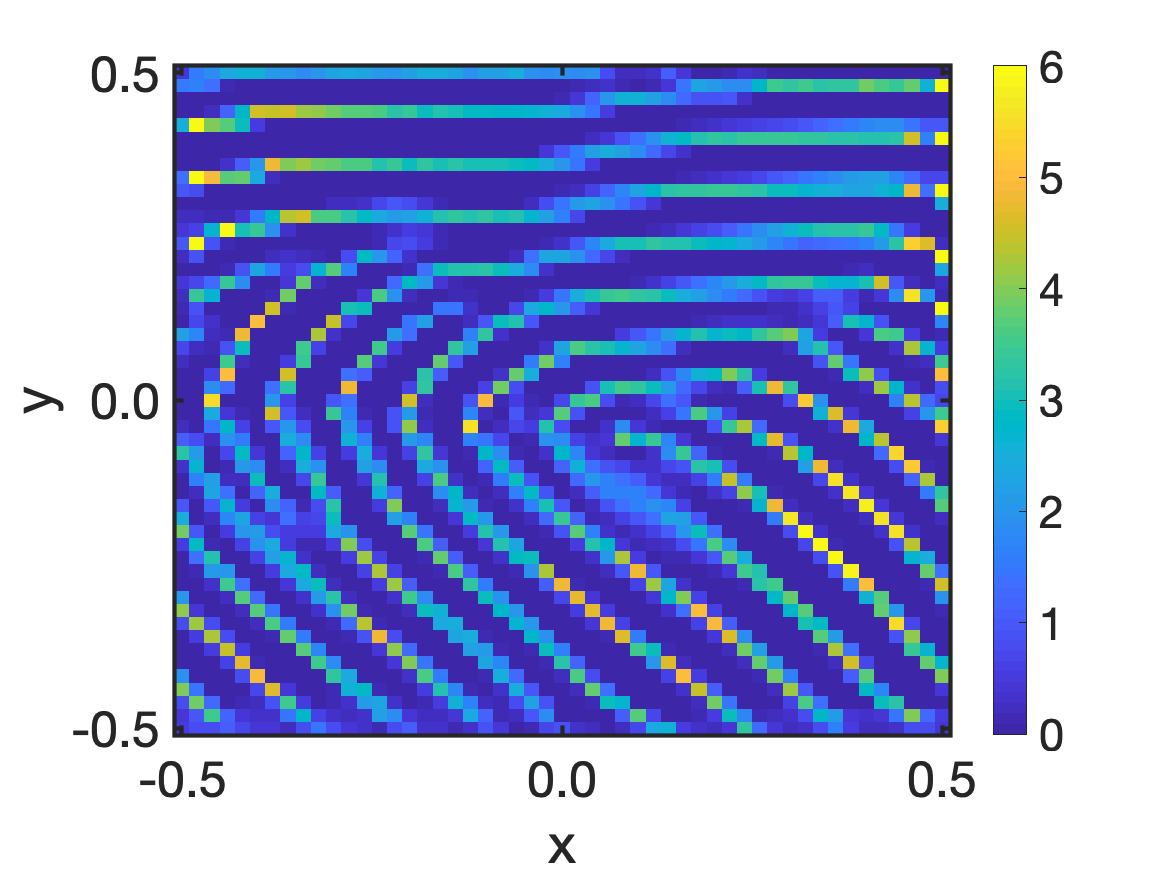}}\\
		\caption{Numerical solution to the anisotropic interaction equation \eqref{eq:macroscopiceqnonlin} after $n$ iterations for different $n$, obtained with the numerical scheme \eqref{eq:numericalscheme} on a grid of size 50 in each spatial direction with diffusion coefficient $\delta=10^{-10}$  for the spatially inhomogeneous tensor field of part of a fingerprint and uniformly distributed initial data on the computational domain $[-0.5,0.5]^2$.}\label{fig:partfingerevolution}
	\end{figure}

	In Figure \ref{fig:partfingerevolution}, we consider the tensor field in Figure \subref*{fig:spart} of part of a fingerprint, and show the numerical solution at different iterations of the numerical scheme \eqref{eq:numericalscheme} on a grid of size 50 in each spatial direction for the diffusion coefficient $\delta=10^{-10}$  and uniformly distributed initial data on the computational domain $[-0.5,0.5]^2$. Note that the resulting numerical solution is close to being stationary.

\begin{figure}[htbp]
		\subfloat[$\delta= 10^{-10}$]{\includegraphics[width=0.32\textwidth]{Meanfield_Partfinger_SpreadInitial_Diffusion1e10}}
		\subfloat[$\delta= 5\cdot 10^{-10}$]{\includegraphics[width=0.32\textwidth]{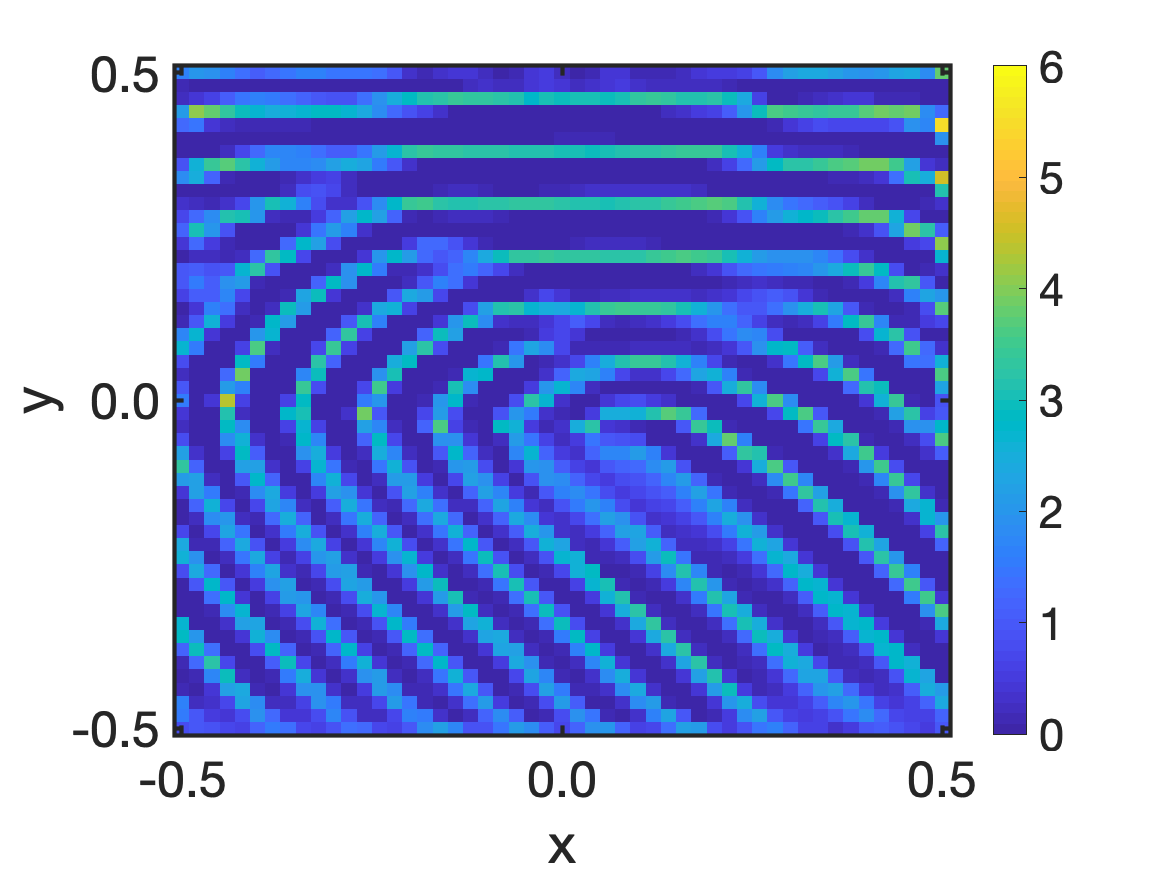}}
		\subfloat[$\delta= 10^{-9}$]{\includegraphics[width=0.32\textwidth]{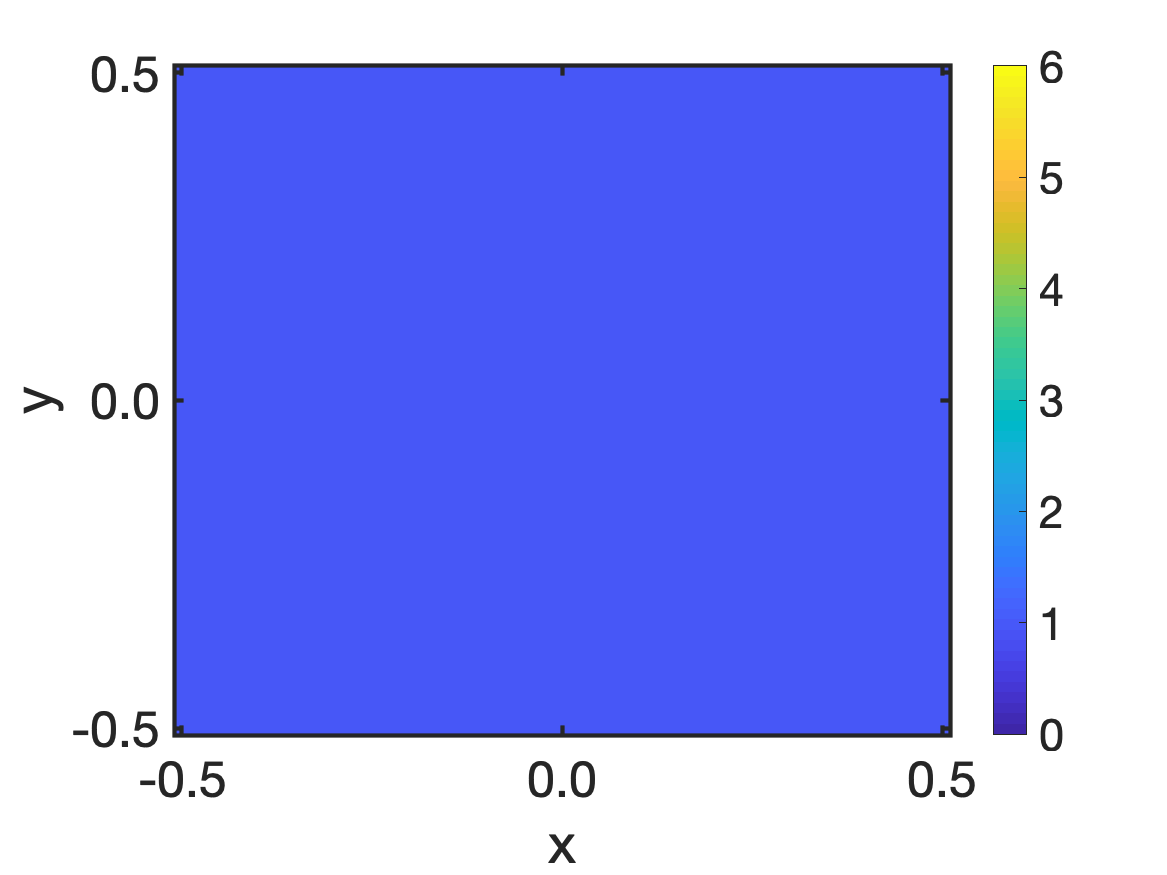}}
		\caption{Stationary solution to the anisotropic interaction equation \eqref{eq:macroscopiceqnonlin}, obtained with the numerical scheme \eqref{eq:numericalscheme} on a grid of size 50 in each spatial direction for different values of the diffusion coefficient $\delta$ for a given spatially inhomogeneous tensor field and uniformly distributed initial data  on the computational domain $[-0.5,0.5]^2$.}\label{fig:spatiallyINhomdiffusionuniformdistributed}
	\end{figure}

\begin{figure}[htbp]
		\subfloat[$\eta= 0.6$]{\includegraphics[width=0.32\textwidth]{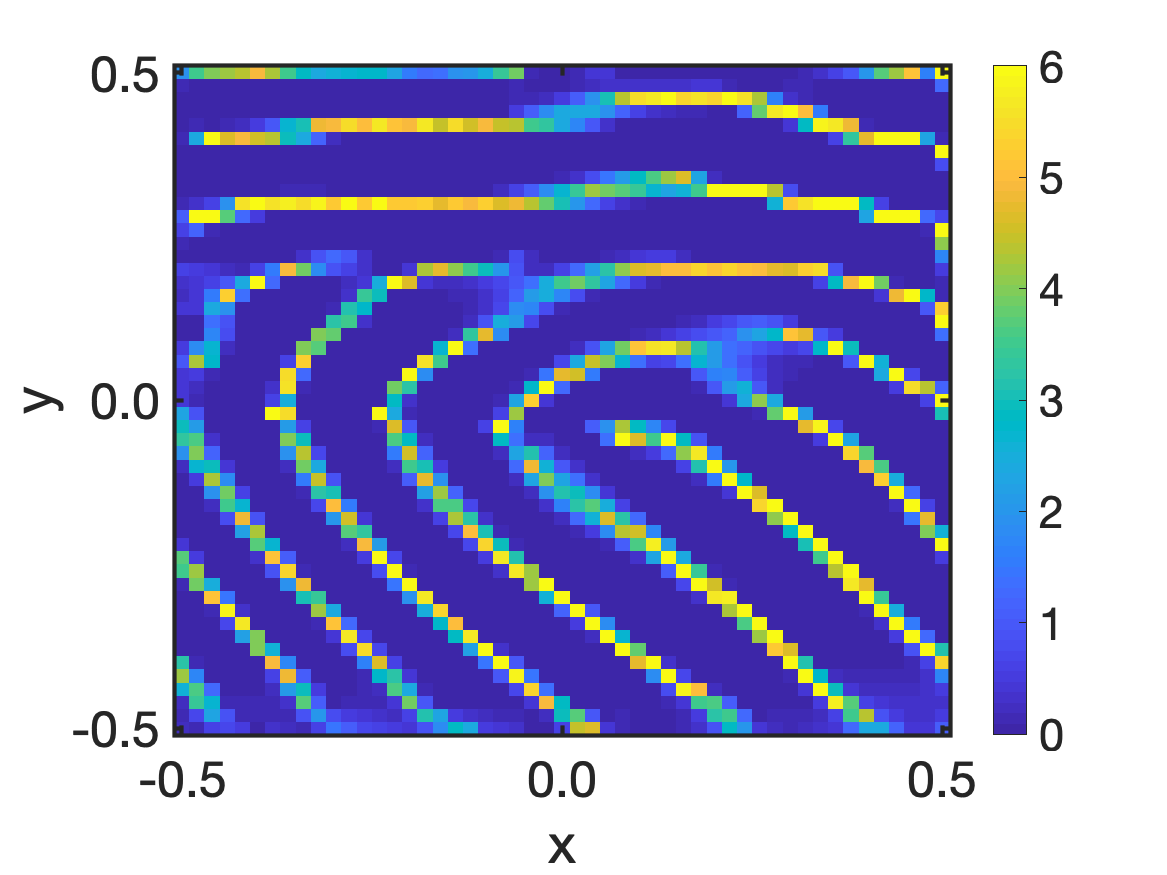}}
		\subfloat[$\eta=0.8$]{\includegraphics[width=0.32\textwidth]{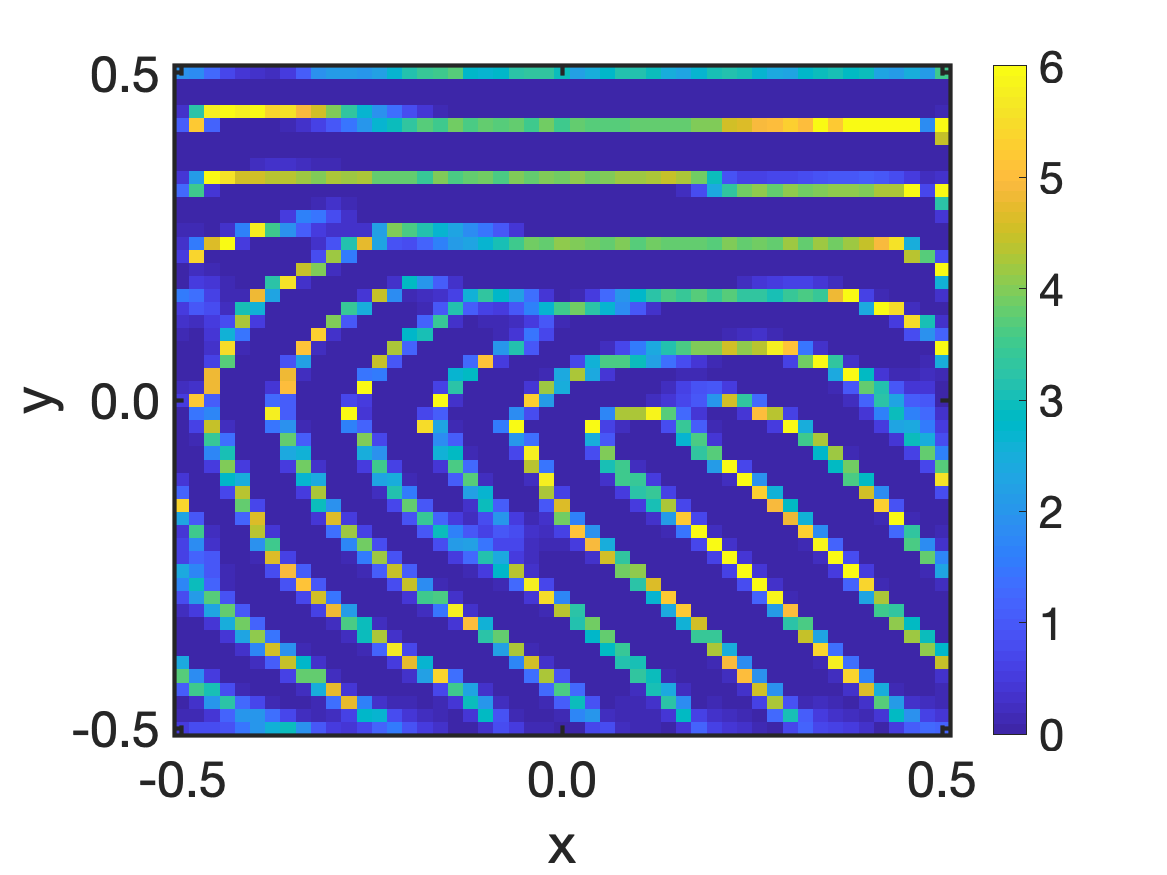}}
		\subfloat[$\eta= 1.2$]{\includegraphics[width=0.32\textwidth]{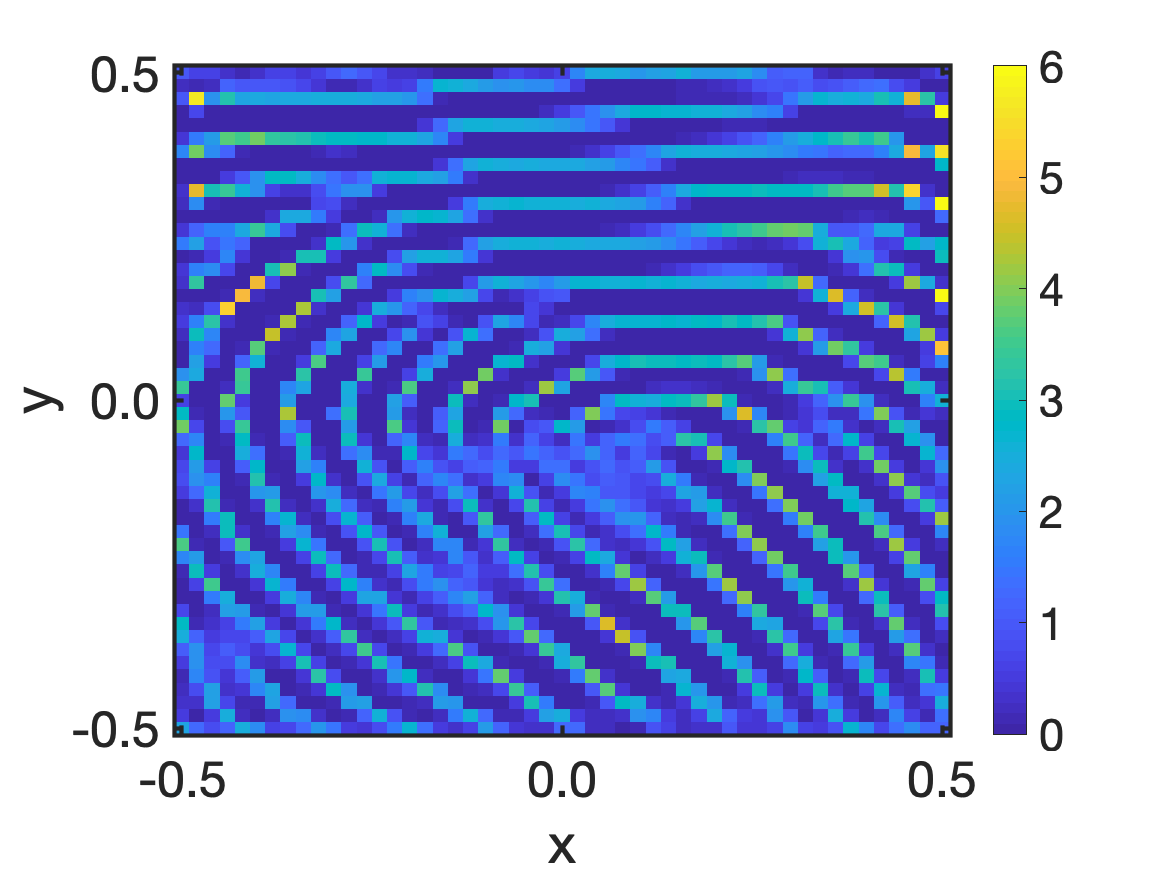}}
		\caption{Stationary solution to the anisotropic interaction equation \eqref{eq:macroscopiceqnonlin}, obtained with the numerical scheme \eqref{eq:numericalscheme} on a grid of size 50 in each spatial direction,  diffusion coefficient $\delta=10^{-10}$ and different force rescalings $\eta$ for a given spatially inhomogeneous tensor field and uniformly distributed initial data  on the computational domain $[-0.5,0.5]^2$.}\label{fig:spatiallyInhomdiffusionuniformdistributed}
	\end{figure}

	Similarly as in Figure \ref{fig:spatiallyhomdiffusionuniformdistributed} for spatially homogeneous tensor fields, we show the stationary solution for different diffusion coefficients $\delta$ in Figure \ref{fig:spatiallyINhomdiffusionuniformdistributed}. For the numerical results in Figure \ref{fig:spatiallyINhomdiffusionuniformdistributed}, the spatially inhomogeneous tensor field in Figure \subref*{fig:spart} and a grid of size 50 in each spatial direction are considered. As $\delta$ increases, the line patterns become wider, provided the diffusion coefficient $\delta$ is below a certain threshold. If $\delta>0$ is above this threshold, e.g.\ for $\delta=10^{-9}$, the uniform distribution is obtained as stationary solution. Note that this threshold is smaller than the one in Figure \ref{fig:spatiallyhomdiffusionuniformdistributed} for spatially homogeneous tensor fields. 
	
	Motivated by the simulation results in \cite{During2017}, we consider different rescalings of the forces in Figure \ref{fig:spatiallyInhomdiffusionuniformdistributed} to vary the distances between the fingerprint lines, i.e.\ we consider $F(\eta d(x,y),T(x))$ where $\eta>0$ is the rescaling factor. As before, we consider the diffusion coefficient $\delta=10^{-10}$ on a grid of size 50 in each spatial direction and uniformly distributed initial data on $[-0.5,0.5]^2$. For $\eta=1$ we recover the same stationary solution  as in Figure \subref*{fig:stationarypart}, while the distances between the fingerprint lines become larger for $\eta\in(0,1)$ and smaller for $\eta>1$. Note that the resulting patterns for the mean-field model  \eqref{eq:macroscopiceqnonlin} are better for $\eta>1$ than for the associated particle model, see \cite[Figure 24]{During2017}, since only dotted lines are possible for particle simulations with $N=2400$ and higher particle numbers result in very long simulation times.

	\section*{Acknowledgments} 	JAC was partially supported by the EPSRC through grant number EP/P031587/1.
	BD has been supported by the Leverhulme Trust research project grant `Novel discretizations for higher-order nonlinear PDE' (RPG-2015-69). LMK was supported by the EPSRC grant Nr.
	EP/L016516/1, the German Academic Scholarship Foundation (Studienstiftung des Deutschen Volkes) and the Cantab Capital Institute for the Mathematics of Information. CBS acknowledges support from the Leverhulme Trust (Breaking the non-convexity barrier, and Unveiling the Invisible), the Philip Leverhulme Prize, the EPSRC grant Nr. EP/M00483X/1, the EPSRC Centre Nr. EP/N014588/1, the European Union Horizon 2020 research and innovation programmes under the Marie Skodowska-Curie grant agreement No. 777826 NoMADS and No. 691070 CHiPS,  the Cantab Capital Institute for the Mathematics of Information  and the Alan Turing Institute.
	
	\bibliographystyle{plain}
	\bibliography{references}
	
\end{document}